\documentclass[numbook]{svjour3}
\usepackage{amsmath,amssymb}
\smartqed

\title{Newton's method on Gra{\ss}mann manifolds}

\newcommand{\C}{\mathbb{C}}
\newcommand{\R}{\mathbb{R}}

\newcommand{\N}{\mathbb{N}}
\newcommand{\half}{\textstyle{\frac{1}{2}}}

\DeclareMathOperator{\colspan}{colspan}
\DeclareMathOperator{\ad}{ad}
\DeclareMathOperator{\T}{T}
\DeclareMathOperator{\Cay}{Cay}
\DeclareMathOperator{\DD}{D}
\DeclareMathOperator{\Gr}{Gr}
\DeclareMathOperator{\LG}{LG}
\DeclareMathOperator{\Hess}{Hess}
\DeclareMathOperator{\e}{e}
\DeclareMathOperator{\id}{id}
\DeclareMathOperator{\im}{im}
\DeclareMathOperator{\tr}{tr}
\DeclareMathOperator{\Sym}{Sym}
\DeclareMathOperator{\GL}{GL}
\DeclareMathOperator{\SO}{SO}
\DeclareMathOperator{\OSp}{OSp}
\DeclareMathOperator{\OO}{O}

\DeclareMathOperator{\U}{U}
\DeclareMathOperator{\so}{\mathfrak{so}}
\DeclareMathOperator{\osp}{\mathfrak{osp}}
\DeclareMathOperator{\vecc}{vec}
\DeclareMathOperator{\dist}{dist}
\DeclareMathOperator{\sinc}{sinc}

\begin{document}
\title{Newton's method on Gra{\ss}mann manifolds\thanks{The first
    author was partially supported by a grant from BMBF within the
    FprofUnd programme. The second and third author are with National
    ICT Australia Limited which is funded by the Australian
    Government's Department of Communications, Information Technology
    and the Arts and the Australian Research Council through
    \emph{Backing Australia's Ability} and the ICT Research Centre of
    Excellence Program.} 
}

\author{Uwe Helmke \and Knut H{\"u}per \and Jochen
  Trumpf\footnote{corresponding author}
}

\institute{U. Helmke \at
  Mathematisches Institut\\
  Universit{\"a}t W{\"u}rzburg\\
  Am Hubland\\
  97074 W{\"u}rzburg\\
  Germany\\
  \email{helmke@mathematik.uni-wuerzburg.de} 
  \and
  K. H{\"u}per \at
  National ICT Australia Limited\\
  Canberra Research Laboratory\\
  Locked Bag 8001\\
  Canberra ACT 2601\\
  Australia\\ 
  \email{knut.hueper@nicta.com.au}\\
  \emph{and}\\
  Department of Information Engineering\\
  The Australian National University
  \and
  J. Trumpf \at
  Department of Information Engineering\\
  The Australian National University\\
  Canberra ACT 0200\\ 
  Australia\\
  Tel.: +61-2-61258677\\
  Fax: +61-2-61258660\\
  \email{jochen.trumpf@anu.edu.au}\\
  \emph{and}\\
  National ICT Australia Limited\\
  Canberra Research Laboratory 
}

\date{Received: date / Accepted: date}

\maketitle

\begin{abstract}
A general class of Newton algorithms on Gra{\ss}mann and Lagrange--Gra{\ss}mann
manifolds is introduced, that depends on an arbitrary pair of local
coordinates. Local quadratic convergence of the algorithm is shown under
a suitable condition on the choice of coordinate systems. Our result extends
and unifies previous convergence results for Newton's method on a manifold. 
Using special choices of the coordinates, new numerical algorithms are derived
for principal component analysis and invariant subspace computations
with improved computational complexity properties.
\keywords{Newton's method \and 
  Gra{\ss}mann and Lagrange Gra{\ss}mann manifold \and 
  smooth parametrizations \and 
  Riemannian metric
}
\subclass{MSC 49M15 \and MSC 53B20 \and MSC 65F15 \and MSC 15A18} 
\end{abstract}

\section{Introduction}
Riemannian optimization is a relatively recent approach towards constrained
optimization that uses full information on the 
underlying geometry of the constraint set in order to set up the optimization algorithms. 
The method is particularly useful if the basic ingredients from differential
geometry, such as the Levi-Civita connection and geodesics are explicitly available. 
This happens in many application problems arising in
signal processing and numerical linear algebra, where optimization naturally 
takes place on homogeneous spaces, such as e.g. Stiefel or Gra{\ss}mann manifolds. In this paper, we 
describe a new class of
Newton algorithms on Gra{\ss}mann manifolds and study applications to 
eigenvalue
and invariant subspace computations.

The idea of using differential geometric methods to construct 
gradient descent algorithms for constrained optimization on smooth manifolds 
is of course not new and we refer to the textbooks
\cite{lue:84,helmke,udriste1} for further information. Such gradient
algorithms use first 
order derivative information on the function and thus can be described in a
rather straightforward way. In contrast, Newton's method on a manifold
requires second order information on the function, using an affine connection 
in order 
to define the Hessian. This can be done in several different ways, thus 
leading to a variety of possible implementations of the Newton algorithm. 

In D. Gabay's work \cite{gabay:82}, the intrinsic Newton method on a
Riemannian manifold is defined via the Levi-Civita connection, taking
iteration steps along associated geodesics. More generally, M. Shub
\cite{shub:86} proposed a Newton method to compute a zero of a smooth
vector field on a smooth manifold endowed with an affine connection. His algorithm is defined for arbitrary families of smooth projections
$\pi_{p}:T_pM\rightarrow M, p\in
  M,$ from the tangent 
bundle which have derivative
equal to the identity at the base point. Therefore it is more general than
Gabay's method and can be employed on arbitrary manifolds, without having to
specify a Riemannian metric. In the case of a gradient vector field on a
Riemannian manifold endowed with the Levi-Civita connection, Shub's algorithm
coincides with Gabay's, when $\{\pi_{p}\}_{p\in M}$ are the Riemannian normal
coordinates.

In the PhD theses of St. Smith and R. Mahony \cite{smi:94a,maho:94a},
see also \cite{ede:98a},  
the Newton method along geodesics of Gabay \cite{gabay:82} was
rediscovered. However, the convergence proofs developed in these 
papers do not apply to the more general situation studied by Shub, 
except for the special case of Rayleigh quotient optimization on the 
unit sphere.  In his recent PhD Thesis,
P.-A. Absil \cite{absil:03a}, see also \cite{absil:04a}, further
discusses the Newton method along 
geodesics and derives a cubic convergence result in a special case.
Moreover, variants with different projections were proposed, too.
There are many more, recent publications discussing aspects of Newton
methods on Riemannian manifolds. We want to specifically mention the
paper by Adler et al. \cite{dedieu} which is similar in spirit to this
paper in so far as it provides explicit formulas for parametrizations and Newton
algorithms on $(\SO_{3})^{N}$. 

In this paper, we propose a general approach to Newton's method on both
Gra{\ss}mann and Lagrange Gra{\ss}mann manifolds that incorporates
the previous ones as special  
cases, but allows also for implementations with improved computational
complexity. We do so by replacing the family of smooth projections by an
arbitrary pair of local coordinates $\mu_p,\nu_p$ with equal derivatives
$D\mu_p(0)=D\nu_p(0)$. Although this generalization might look minor at first
sight, it is actually crucial to achieve better performance. Following 
\cite{huep:05a} and extending the known local quadratic convergence result for
the intrinsic Riemannian Newton method, we prove local quadratic convergence
of the generalized Newton algorithm. The Newton method on the Lagrange
Gra{\ss}mannian has not been considered before, but has important 
applications in control (e.g. to algebraic Riccati equations in linear 
quadratic control). 

The paper is structured as follows.
In order to enhance the readability of the paper for non-experts, we begin
with a brief summary of the basic differential geometry of the classical 
Gra{\ss}mann
manifold and the Lagrange Gra{\ss}mannian, respectively, deriving explicit formulas for (projections onto) 
tangent spaces, normal spaces, gradients, Hessians, and geodesics. 
We then compute the Riemannian normal coordinates of the two types of
Gra{\ss}mannians. Using approximations of the exponential map via e.g. Pad{\'e}
approximants or the $QR$ factorization, then leads to alternative coordinate
systems and resulting simplified implementations of the Newton algorithm.
By generalizing the construction of Shub, we introduce the Newton algorithm
via a pull back/push forward scheme defined by an arbitrary pair of local
coordinates for the Gra{\ss}mannians. This leads to a rich family of 
intrinsically
defined Newton methods that have potential for considerable
computational advantages compared with the 
previously known algorithms. In fact, instead of relying upon the use of
Riemannian normal coordinates, that are difficult to compute with, we 
advocate to use the much more easily computable local coordinates 
via the $QR$-factorization.

For example, in Edelman et al. \cite{ede:98a} the steps of the Newton algorithm on the
classical Gra{\ss}mannian are defined in the ambient Euclidean space of the
associated Stiefel manifold. This leads them to solving sequences of Sylvester
equations in higher dimensional matrix spaces than necessary. In contrast, our 
algorithms works with the minimal number of parameters, given by the dimension
of the Gra{\ss}mannian. Moreover, our algorithms do not require the
iterative calculation of matrix exponentials, but only involve finite step
iterations using efficient $QR$-computations.

Finally, we apply these techniques to eigenspace computations. By applying our
Newton scheme to the Rayleigh quotient
function on the Gra{\ss}mann (and Lagrange Gra{\ss}mann) manifold, we obtain
a new class of iterative algorithms for principal component analysis with
improved computational complexity. For eigenspace computations of arbitrary,
not necessarily symmetric, matrices we derive an apparently new class of
Newton algorithms, that requires the repeated computations of solutions 
to nested
Sylvester type equations.

\section{Riemannian geometry of  the Gra{\ss}mann manifold}

In this section we describe the basics for the Riemannian geometry
of Gra{\ss}mann manifolds, i.e. tangent and normal spaces, Riemannian
metrics and geodesics. We focus on the real Gra{\ss}\-mannian; the
results carry through mutatis mutandis for complex Gra{\ss}\-mannians,
too.

Recall, that the Gra{\ss}mann manifold $\Gr_{m,n}$ is
defined as the set of $m$-dimensional $\R$-linear subspaces of
$\R^n$. It is a smooth, compact manifold of dimension $m(n-m)$ and
provides a natural generalization of the familiar projective
spaces. Let denote
\begin{equation}
  \label{eq:5}
\OO_{n}:=\{X\in\R^{n\times n}|X^{\top}X=I\}.
\end{equation}
and
\begin{equation}
  \label{eq:3}
  \SO_{n}:=\{X\in\OO_{n}|\det X=1\}
\end{equation}
The Gra{\ss}mann manifold can also be viewed in an equivalent way as a
homogeneous space $\mathrm{SO}_{n}(\R)/H$, cf. e.g. \cite{helmke} and
see below for a definition of $H$, for the transitive
$\SO_{n}$--action
\begin{equation}\label{eq:10}
  \begin{split}
    \sigma:\SO_{n}\times\Gr_{m,n}&\to
    \Gr_{m,n},\\
    (T,\mathsf{V})&\mapsto T\mathsf{V}.
  \end{split}
\end{equation}
Let
\begin{equation}\label{eq:2}
  \mathsf{V_{0}}=\colspan
  \begin{bmatrix}
    I_{m} \\ 0
  \end{bmatrix} \in \Gr_{m,n}
\end{equation}
denote the standard $m$-dimensional subspace of $\R^n$ that is
spanned by the first $m$ standard basis vectors of $\R^n$. Then
the stabilizer subgroup $H:=\mathrm{Stab}(\mathsf{V}_{0})$ of
$\mathsf{V}_{0}$ is given by
\begin{equation}\label{eq:H}
  H =\left\{
    \begin{bmatrix}
      U & 0 \\
      0 & V
    \end{bmatrix}\in \SO_{n}\,\left|\,
    {U\in \mathrm{O}_{m}},\ {V\in \mathrm{O}_{n-m}}
  \right.\right\},
\end{equation}
i.e. by the compact Lie subgroup of $\SO_{n}$
consisting of all block diagonal orthogonal matrices. The map
\begin{equation}
  \label{eq:6}
\SO_{n}/H \to \Gr_{m,n},\ {\Theta}H \mapsto {\Theta}\mathsf{V}_{0}  
\end{equation}
then defines a diffeomorphism of the Gra{\ss}mann manifold with the
homogeneous space $\SO_{n}/H$. See Edelman et al. \cite{ede:98a},
Absil \cite{absil:03a} and H\"uper and Trumpf \cite{huep:05a} for
further details on Newton's method on
$\Gr_{m,n}$, in a variant that exploits the homogeneous space
structure of the Gra{\ss}mann manifold. Here we develop a different
approach, by identifying $\Gr_{m,n}$ with a set
of self-adjoint projection operators. 

Thus we define the {\em Gra{\ss}mannian} as
\begin{equation}
  \label{eq:7}
  \Gr_{m,n}:= \{P \in \R^{n\times n}\ | \ P^{\top}=P, P^2=P,\tr P=m\},
\end{equation}
the manifold of rank $m$ symmetric projection operators of $\R^n$;
see \cite{helmke} for the con\-struction of a natural bijection
with the Gra{\ss}mann manifold and a proof that it defines a
diffeomeorphism. In the sequel we will describe the Riemannian
geometry directly for the submanifold $\Gr_{m,n}$
of $\R^{n\times n}$. As we will see, this approach has advantages
that simplify both the analysis and design of Newton-based
algorithms for the computation of principal components.

We begin by recalling the following known and basic fact on the
Gra{\ss}mannian; see \cite[Section 2.1]{helmke} for a proof in the
more general context of isospectral manifolds. Let
\begin{equation}
  \label{eq:8}
 \Sym_{n}:=\{S\in \R^{n\times n}\ | \ S^{\top}=S\} 
\end{equation}
and
\begin{equation}
  \label{eq:9}
\so_{n}:=\{\Omega\in \R^{n\times n}\ | \ \Omega^{\top}=-\Omega\} 
\end{equation}
denote the vector spaces of
real symmetric and real skew-symmetric matrices, respectively.
\begin{theorem}
\begin{enumerate}
\item[(a)]The Gra{\ss}mannian $\Gr_{m,n}$ is a smooth,
compact submanifold of $\Sym_{n}$ of dimension $m(n-m)$.
\item[(b)] The tangent space of $\Gr_{m,n}$ at an
element $P\in \Gr_{m,n}$ is given as
\begin{equation}
  \label{eq:11}
  T_P\Gr_{m,n}=\{[P,\Omega]\ | \ \Omega\in \so_{n}\}.
\end{equation}
Here $[P,\Omega]:=P\Omega - \Omega P$ denotes the matrix commutator
(Lie bracket). 
\end{enumerate}
\end{theorem}

Let
\begin{equation}
  \label{eq:12}
  \begin{split}
 \ad_P:\R^{n\times n} &\to \R^{n\times n},\\
\ad_P(X)&:=[P,X]   
  \end{split}
\end{equation}
denote the adjoint representation at $P$. For a
projection operator $P$ it enjoys the following property.
\begin{lemma} 
\label{lem:adP3}
For any $P\in \Gr_{m,n}$, the minimal polynomial of
$\ad_P : \R^{n\times n} \to \R^{n\times n}$ is equal to
$s^3 - s$. Thus $\ad_{P}^3=\ad_{P}$. 
Moreover, 
\begin{equation}
  \label{eq:13}
  \ad_{P}^2X =[P,[P,X]]=X
\end{equation}
holds for all tangent vectors $X \in
T_{P}\Gr_{m,n}$.
\end{lemma}

\begin{proof} 
From $P^2=P$ we get
\begin{equation}
  \label{eq:100}
  \ad_{P}^2X=[P,[P,X]]=P^2X+XP^2-2PXP=PX+XP-2PXP
\end{equation}
and therefore, using $P^2=P$ again
\begin{equation}
  \begin{split}
    \ad_{P}^3X&=P(PX+XP-2PXP)-(PX+XP-2PXP)P\\&=PX-XP\\&=\ad_{P}X  
  \end{split}
\label{eq:101}
\end{equation}
for all $n\times n$--matrices $X$. If $X=[P,\Omega]$ is a tangent
vector, then $\ad_{P}^2X=\ad_{P}^3\Omega=\ad_{P}\Omega=X$. The result
follows. 
\qed\end{proof}

We use this result to describe the normal bundle of
$\Gr_{m,n}$. In the sequel, we will always endow
$\Sym_{n}$ with the Frobenius inner product, defined by
\begin{equation}
  \label{eq:98}
 \langle X,Y\rangle := \tr(XY) 
\end{equation}
for all $X,Y\in \Sym_{n}$. Since the tangent space
$T_{P}\Gr_{m,n}\subset \Sym_{n}$ is a subset
of $\Sym_{n}$ (using the usual identification of $T_{P}\Sym_{n}$ with
$\Sym_{n}$), we can define the normal space at $P$ to be 
the vector space
\begin{equation}
  \label{eq:99}
  N_{P}\Gr_{m,n}=\!(T_{P}\Gr_{m,n})^{\perp} := \{X\in \Sym_{n}|\tr(XY)=0
\ \text{for all}\ Y\in T_{P}\Gr_{m,n}\}.
\end{equation}

\begin{proposition} 
\label{prop:normalGR}
Let $P\in \Gr_{m,n}$ be arbitrary.
\begin{enumerate}
\item 
The normal subspace in $\Sym_{n}$ is given as
\begin{equation}
  \label{eq:48}
  N_{P}\Gr_{m,n}=\{X-\ad_{P}^2X \ | \ X\in \Sym_{n}\}.
\end{equation}
\item
The linear map
\begin{equation}
  \label{eq:49}
  \pi: \Sym_{n}\to \Sym_{n},\ X\mapsto \ad_{P}^{2}X=[P,[P,X]]
\end{equation}
is the self-adjoint projection operator onto
$T_{P}\Gr_{m,n}$ with kernel $N_{P}\Gr_{m,n}$.
\end{enumerate}
\end{proposition}

\begin{proof} For any tangent vector $[P,\Omega]\in
T_{P}\Gr_{m,n}$, where $\Omega^{\top}=-\Omega$, and
any $X=X^{\top}$, we have
\begin{equation}
\begin{split}
\tr([P,\Omega](X-\ad_{P}^2X))&=\tr(([X,P]-[\ad_{P}^2X,P])\Omega)\\
&=\tr(([X,P]+\ad_{P}^3X)\Omega)\\
 &=\tr((\ad_{P}^3X - \ad_{P}X)\Omega)\\
&= 0,
\end{split}
\end{equation}
since $\ad_{P}^3=\ad_{P}.$ 
Therefore,
$T_{P}\Gr_{m,n}$ and $\{X-\ad_{P}^2X \ | \ X\in \Sym_{n}\}$ are
orthogonal subspaces of $\Sym_{n}$ 
with respect to the Frobenius inner product. Their sum also spans
$\Sym_{n}$, as otherwise there exists a nontrivial $S \in \Sym_{n}$
that is orthogonal to both spaces; but then for all $\Omega \in
\so_{n}$
\begin{equation}
  \label{eq:4}
  \tr(S[P,\Omega])=\tr([S,P]\Omega)=0\quad \Longrightarrow\quad 
[S,P]=0,
\end{equation}
and for all $X\in \Sym_{n}$, using \eqref{eq:4}
\begin{equation}
\tr(S(X-\ad_{P}^2X))=\tr(SX-[S,P][P,X])=\tr(SX)=0
\end{equation}
which implies $S=0$, a contradiction. Thus
the two spaces define an orthogonal sum decomposition of
$\Sym_{n}$ and therefore $\{X-\ad_{P}^2X  |  X\in
\Sym_{n}\}$ must be the normal space. This completes the
proof for the first claim. 

Since $\pi=\ad_{P}^2$, we have
\begin{equation}
  \label{eq:16}
  \pi^{2}=\ad_{P}^4=\ad_{P}^2=\pi
\end{equation}
because $\ad_{P}^3=\ad_{P}$.
Moreover, by
definition of $\pi$ we have 
$\im \pi \subset T_{P}\Gr_{m,n}$, cf. \eqref{eq:11}, and for any $X\in
T_{P}\Gr_{m,n}$ we have by Lemma~\ref{lem:adP3} that $\pi(X)=X$.
Therefore
\begin{equation}
  \label{eq:17}
  \im \pi =
T_{P}\Gr_{m,n}.
\end{equation}
For any $X-\ad_{P}^{2}X \in N_{P}\Gr_{m,n}$ 
we have
\begin{equation}
  \label{eq:18}
  \pi(X-\ad_{P}^2X)=\ad_{P}^2X -
\ad_{P}^4X=0,
\end{equation}
by~\eqref{eq:16}. Since $N_{P}\Gr_{m,n}$ is the
orthogonal complement to the tangent space in $\Sym_{n}$, a straight forward dimension
argument yields $\ker \pi = N_{P}\Gr_{m,n}$. Finally, using the
Frobenius inner product on $\Sym_{n}$, we have for all $X_{1},X_{2}\in\Sym_{n}$
\begin{equation}
  \begin{split}
 \langle \pi(X_{1}),X_{2}\rangle&=\tr((\ad_{P}^2 X_{1})X_{2})\\
&=\tr([P,[P,X_{1}]]X_{2})\\
&=\tr([P,[P,X_{2}]]X_{1})\\
&=\langle X_{1},\pi(X_{2})\rangle.  
  \end{split}
\end{equation}
Thus $\pi$ is self-adjoint and the result follows. 
\qed\end{proof}

A formula for $\pi$ in the language of linear maps has already been
given in \cite[Section 4.2]{machadosalavessa:85}.

There are at least two natural Riemannian metrics defined on the
Gra{\ss}mannian $\Gr_{m,n}$, the induced {\em Euclidean
metric} and the {\em normal metric}, cf. e.g. \cite{helmke} or
\cite{maho:94a}.  

The Euclidean Riemannian
metric on $\Gr_{m,n}$ is defined by the Frobenius
inner product on the tangent spaces
\begin{equation}
  \label{eq:97}
  \langle X,Y \rangle := \tr(XY)
\end{equation}
for all $X,Y\in T_{P}\Gr_{m,n}$ which is induced by the embedding space
$\Sym_{n}$. 

The normal Riemannian metric has a
somewhat more complicated definition. Consider the surjective
linear map
\begin{equation}
  \label{eq:19}
  \begin{split}
    \ad_P: \so_{n} &\to T_{P}\Gr_{m,n},\\ \Omega&\mapsto [P,\Omega]
  \end{split}
\end{equation}
with kernel
\begin{equation}
  \label{eq:102}
  \ker\ad_P = \{\Omega \in \so_{n} \ | \ P\Omega = \Omega P\}.
\end{equation}
We regard $\so_{n}$ as an inner product space, endowed
with the Frobenius inner product $\langle \Omega_1, \Omega_2
\rangle = \tr(\Omega_{1}^{\top}\Omega_2) =
-\tr(\Omega_{1}\Omega_{2})$.
Then $\ad_P$ induces an
isomorphism of vector spaces
\begin{equation}
  \label{eq:55}
  \widehat{\ad}_P: (\ker \ad_P)^{\perp} \to T_{P}\Gr_{m,n}
\end{equation}
and therefore induces an isometry of inner product spaces, by
defining an inner product on $T_{P}\Gr_{m,n}$ via
\begin{equation}
  \label{eq:56}
  \langle\langle X,Y \rangle\rangle _{P} := -\tr(\widehat{\ad}_{P}^{-1}(X)
\widehat{\ad}_P^{-1}(Y)).
\end{equation}
Note, that this inner product on $T_{P}\Gr_{m,n}$,
called the normal Riemannian metric, might vary with the basepoint
$P$. Luckily, the situation is better than one would expect, as
Proposition~\ref{prop:equalmetrics} below shows.

But first we will show that the operator $\ad_{P}^{2}$,
$P\in\Gr_{m,n}$, is equally well behaved on $\so_{n}$ as it is on
$\Sym_{n}$, cf. Proposition~\ref{prop:normalGR}.

\begin{proposition}
  \label{prop:adP2son}
  Let $P\in\Gr_{m,n}$ be arbitrary. The linear map
  \begin{equation}
    \label{eq:90}
    \ad_{P}^2: \so_{n}\to \so_{n}, \ \Omega\mapsto[P,[P,\Omega]]
  \end{equation}
  is the self-adjoint projection operator onto $(\ker\ad_{P})^{\perp}$
  along $\ker\ad_{P}$.
\end{proposition}

\begin{proof}
  Let $\Omega\in\so_{n}$ be arbitrary. 
  By Lemma~\ref{lem:adP3} we know that 
  $X:=\ad_P(\Omega)=\ad_{P}(\ad_{P}^2\Omega)$. But since 
  for all $\Omega_{1}\in\ker\ad_P\subset\so_{n}$ 
  \begin{equation}
    \label{eq:20}
    \tr(\ad_{P}^2(\Omega)\Omega_{1}^{\top})=\tr(\ad_P(\Omega)\ad_{P}(\Omega_{1}))=0, 
  \end{equation}
  we conclude
  $\ad_{P}^2\Omega\in (\ker\ad_P)^{\perp}\subset\so_{n}$. 
  Now let $\Omega\in(\ker\ad_{P})^{\perp}$. Then $\ad_{P}^2\Omega\in
  (\ker\ad_P)^{\perp}$ and hence $\Omega-\ad_{P}^2\Omega\in
  (\ker\ad_P)^{\perp}$. By Lemma~\ref{lem:adP3} 
  $\ad_{P}(\Omega-\ad_{P}^2\Omega)=\ad_{P}\Omega-\ad_{P}^{3}\Omega=0$
  and hence $\Omega-\ad_{P}^2\Omega\in\ker\ad_{P}$. It follows 
  $\Omega-\ad_{P}^2\Omega=0$ and thus $\ad_{P}^{2}\Omega=\Omega$.

  We have shown $\im\ad_{P}^{2}\subset(\ker\ad_P)^{\perp}$ and that
  the restriction of $\ad_{P}^{2}$ to $(\ker\ad_P)^{\perp}$ is the
  identity. It remains to show that $\ker\ad_{P}^{2}=\ker\ad_{P}$, but
  this follows readily from Lemma~\ref{lem:adP3}.
\qed\end{proof}

\begin{proposition}\label{prop:equalmetrics}The Euclidean and normal Riemannian metrics on
the Gra{\ss}mannian $\Gr_{m,n}$ coincide, i.e. for all
$P\in \Gr_{m,n}$ and for all $X,Y \in
T_{P}\Gr_{m,n}$ we have
\begin{equation}
  \label{eq:57}
  \tr(X^{\top}Y)=-\tr\left(\widehat{\ad}_{P}^{-1}(X)\ \widehat{\ad}_P^{-1}(Y)\right).
\end{equation}
\end{proposition}

\begin{proof} 
Choose
\begin{equation}
  \label{eq:87}
  \Omega_1, \Omega_2 \in (\ker\ad_P)^{\perp}\;\text{with}\;
X=[P,\Omega_1]\;\text{and}\; Y=[P,\Omega_2].
\end{equation}
Then
\begin{equation}
  \label{eq:107}
  -\tr(\widehat{\ad}_{P}^{-1}(X)\widehat{\ad}_P^{-1}(Y))=
\tr(\Omega_{1}^{\top}\Omega_{2}).
\end{equation}
On the other hand
\begin{equation}
  \begin{split}
  \tr(X^{\top}Y)&=\tr([P,\Omega_{1}][P,\Omega_{2}])\\
&=\tr([P,[P,\Omega_{1}]]^{\top}\Omega_{2}).  
  \end{split}
\end{equation}
Now by Proposition~\ref{prop:adP2son} we know that 
$\ad_{P}^{2}\Omega_{1}=\Omega_{1}$
and this implies
\begin{equation}
  \label{eq:21}
  \tr(X^{\top}Y)=\tr(\Omega_{1}^{\top}\Omega_{2}),
\end{equation}
as claimed. 
\qed\end{proof}

Since these two Riemannian metrics on the Gra{\ss}mannian coincide,
they also define the same geodesics. Thus, in the sequel, we focus
on the Euclidean metric. Note, that the above result is not true
for arbitrary flag manifolds and in fact, the geodesics are then
different for the two metrics. The following result characterizes
the geodesics on $\Gr_{m,n}$.
\begin{theorem} 
\label{theo:GeoGrass}
The geodesics of $\Gr_{m,n}$ are
exactly the solutions of the second order differential equation
\begin{equation}\label{eq:geodesic1}
\ddot{P} + [\dot{P},[\dot{P},P]] = 0.
\end{equation}
The unique geodesic $P(t)$ with initial conditions $P(0)=P_0 \in
\Gr_{m,n}$, $\dot{P}(0)=\dot{P}_0 \in
T_{P_0}\Gr_{m,n}$ is given by
\begin{equation}\label{eq:geodesic2}
P(t) = \e^{t[\dot{P}_0,P_0]}P_0\e^{-t[\dot{P}_0,P_0]}.
\end{equation}
\end{theorem}
\begin{proof} The geodesics of $\Gr_{m,n}$ for the
Euclidean metric are characterized as the curves $P(t)\in
\Gr_{m,n}$, such that $\ddot{P}(t)$ is a normal
vector for all $t\in \R$. This condition is equivalent to the
existence of $S(t)=S(t)^{\top}$ with
\begin{equation}\label{eq:normal1}
\ddot{P}=S-\ad_P^{2}S.
\end{equation}
Since by Lemma~\ref{lem:adP3} $\ad_P^{3} = \ad_P$ this
implies
\begin{equation}
  \label{eq:22}
  \ad_P\ddot{P}=[P,\ddot{P}]=0.
\end{equation}
Moreover, any curve $P(t)\in \Gr_{m,n}$ satisfies the
identity
\begin{equation}\label{eq:normal2}
\dot{P}=\ad_P^{2}\dot{P},
\end{equation}
as $\ad_{P}^2$ acts as the identity on the tangent space
$T_{P}\Gr_{m,n}$. By differentiating equation
(\ref{eq:normal2}) we obtain
\begin{equation}
  \label{eq:58}
  \begin{split}
  \ddot{P}&=[P,[P,\ddot{P}]]+[P,[\dot{P},\dot{P}]]+[\dot{P},[P,\dot{P}]]\\
&=
-\ad_{\dot{P}}^{2}P + \ad_{P}^{2}\ddot{P}.  
  \end{split}
\end{equation}
Therefore, if $P(t)$ is a geodesic, then $\ad_{P}^2(\ddot{P})=0$ and
\begin{equation}
  \label{eq:59}
  \begin{split}
  \ddot{P}&=-\ad_{\dot{P}}^{2}P + \ad_{P}^{2}\ddot{P}\\
&= -\ad_{\dot{P}}^{2}P
  \end{split}
\end{equation}
and therefore satisfies $\ddot{P} + [\dot{P},[\dot{P},P]]=0$, as
claimed. We now check, that every curve $P(t)$ as in
(\ref{eq:geodesic2}) is a solution to (\ref{eq:geodesic1}). Let
$\Omega :=[\dot{P}_{0},P_{0}]$. Then
\begin{equation}
  \label{eq:60}
  \dot{P}=[\Omega,P], \quad \ddot{P}=[\Omega,[\Omega,P]]
\end{equation}
and thus (\ref{eq:geodesic1})is equivalent to
\begin{equation}\label{eq:geodesic3}
[\Omega,[\Omega,P]]+[[\Omega,P],[[\Omega,P],P]]=0.
\end{equation}
Multiplying by the left and right with $\e^{-t\Omega}$ and
$\e^{t\Omega}$ respectively, we see that (\ref{eq:geodesic3}) is
equivalent to
\begin{equation}
[\Omega,[\Omega,P_0]]+[[\Omega,P_0],[[\Omega,P_0],P_0]]=0.
\end{equation}
Without loss of generality we can assume that
\begin{equation}
  \label{eq:61}
  P_{0}=\begin{bmatrix}
      I_m & 0 \\
      0 &  0
    \end{bmatrix}
\end{equation}
and therefore
\begin{equation}
  \label{eq:62}
  \Omega=\begin{bmatrix}
      0 & Z \\
      -Z^{\top} & 0
    \end{bmatrix}
\end{equation}
with $Z\in\R^{m\times(n-m)}$. Thus
\begin{equation}
  \label{eq:63}
  [\Omega,[\Omega,P_0]]=\begin{bmatrix}
      -2ZZ^{\top} & 0 \\
       0 & 2Z^{\top}Z
    \end{bmatrix}
\end{equation}
and also
\begin{equation}
  \label{eq:64}
  [[\Omega,P_0],[[\Omega,P_0],P_0]]=\begin{bmatrix}
      2ZZ^{\top} & 0 \\
       0 & -2Z^{\top}Z
    \end{bmatrix}.
\end{equation}
This implies (\ref{eq:geodesic3}) and shows that any curve given
by (\ref{eq:geodesic2}) is a solution of (\ref{eq:geodesic1}).
Since any $P_0\in \Gr_{m,n}$ and $[\dot{P}_0,P_0]\in
T_{P}\Gr_{m,n}$ are admissible initial conditions for
(\ref{eq:geodesic1}), and since the resulting initial value problem
has a unique solution (namely \eqref{eq:geodesic2}), this shows that (\ref{eq:geodesic2}) is
exactly the set of {\em all} solutions of (\ref{eq:geodesic1}).
Moreover, for the particular initial point
\begin{equation}
  \label{eq:54}
  P_0=\begin{bmatrix}
      I_m & 0 \\
      0 &  0
    \end{bmatrix}
\end{equation}
one observes that
\begin{equation}
  \label{eq:65}
  [[\Omega,P_0],[[\Omega,P_0],P_0]]=\begin{bmatrix}
      2ZZ^{\top} & 0 \\
       0 & -2Z^{\top}Z
    \end{bmatrix}
\end{equation}
is a normal vector to the Gra{\ss}mannian at $P_0$. Thus, by
invariance of the normal bundle under orthogonal similarity
transformations $P_0 \mapsto \Theta^{\top} P_{0}\Theta,\ \Theta\in
\SO_{n}$, we see that $[\dot{P},[\dot{P},P]]$ is a
normal vector to $T_{P}\Gr_{m,n}$ for all $\dot{P}\in
T_{P}\Gr_{m,n}$. Thus, for any solution $P(t)$ of
(\ref{eq:geodesic1}) also $\ddot{P}=-[\dot{P},[\dot{P},P]]$ is a
normal vector, and hence all solutions of (\ref{eq:geodesic1}) are
geodesics. 
\qed\end{proof}

The above explicit formula for geodesics leads to the following
formula for the geodesic distance between two points on a
Gra{\ss}mannian. We omit the simple proof; see also \cite{absil:03a}
for a slightly different formula which is only valid on an open and
dense subset
of the Gra{\ss}mannian.

\begin{corollary}
  Let $P,Q\in\Gr_{m,n}$. Given any $\Theta\in\SO_{n}$ such that 
 \begin{equation*}
  P=\Theta^{\top}\begin{bmatrix}
      I_m & 0 \\
      0 &  0
    \end{bmatrix}\Theta
\end{equation*} we define
\begin{equation*}
    \begin{bmatrix}Q_{11} & Q_{12} \\ Q_{12}^{\top} & Q_{22}\end{bmatrix}:=
  \Theta Q\Theta^{\top}.
\end{equation*}
Let $1\geq\lambda_{1}\geq\cdots\geq\lambda_{m}\geq 0$ denote the
eigenvalues of $Q_{11}$. The geodesic distance of
$P$ to $Q$ in $\Gr_{m,n}$ is given by
\begin{equation}
  \label{eq:76}
  \dist(P,Q)=\sqrt{2\sum_{i=1}^{m}\arccos^{2}(\sqrt{\lambda_{i}})}\ .
\end{equation}
Alternatively, let $1\geq\mu_{1}\geq\cdots\geq\mu_{n-m}\geq 0$ denote
the eigenvalues of $Q_{22}$. Then
\begin{equation}
  \label{eq:77}
  \dist(P,Q)=\sqrt{2\sum_{i=1}^{n-m}\arcsin^{2}(\sqrt{\mu_{i}})}\ .
\end{equation}
In particular, if $P,Q\in\Gr_{m,n}$ with $Q=YY^{\top},Y^{\top}Y=I_m $, then
\begin{equation}
  \label{eq:76a}
  \frac{1}{2}\dist^2(P,Q)=\tr
    \left(\arccos^{2}((Y^{\top}PY)^{\frac{1}{2}})
    \right)\ .
\end{equation}
\end{corollary}

Note that formula~\eqref{eq:77} is more efficient in the case $2m>n$.
Note also that our formulas imply that the maximal length of a simple closed
geodesic in $\Gr_{m,n}$ is $\sqrt{2m}\cdot\pi$ for $2m\leq n$ and
$\sqrt{2(n-m)}\cdot\pi$ for $2m>n$.

\subsection{Parametrizations and Coordinates for the Gra{\ss}mannian}
\label{sec:param-coord-grassm}
In this section we briefly recall the notion of local parametrization for
smooth manifolds. For further details we refer to \cite{Lang:99}. Let $M$ be
a smooth $n$-dimensional real manifold then for every point $p\in M$ there
exists a smooth map
\begin{equation*}
  \mu_{p}:\R^{n}\longrightarrow M, \ \mu_{p}(0)=p
\end{equation*}
which is a local diffeomorphism around $0\in\R^{n}$. Such a map is called
a \emph{local parametrization around $p$}.

We consider local parametrizations for the
Gra{\ss}mannian via the tangent space, i.e. families of smooth maps
\begin{equation}
  \label{eq:66}
 \mu_P : T_{P}\Gr_{m,n}\to \Gr_{m,n} 
\end{equation}
satisfying
\begin{equation}
  \label{eq:108}
  \mu_{P}(0)=P\;\text{and} \;\DD\mu_{P}(0)=\id.
\end{equation}
We introduce three
different choices of such local parametrizations.

\subsubsection{Riemannian normal coordinates}
\label{sec:riem-norm-coord}
Riemannian normal coordinates are defined through the Riemannian
exponential map (see e.g. \cite{jost2})
\begin{equation}\label{eq:Rinoco}
  \begin{split}
    \mu_{P}^{\exp}=\exp_{P}: T_{P}\Gr_{m,n} &\to \Gr_{m,n},\\
    \exp_{P}(\xi) &= \e^{[\xi,P]}P\e^{-[\xi,P]}.
  \end{split}
\end{equation}

\begin{remark}
  \label{sec:riem-norm-coord-2}
Note that by Theorem~\ref{theo:GeoGrass} the unique geodesic $P(t)$ with
initial conditions $P(0)=P_{0}$ and $\dot{P}(0)=\dot{P}_{0}$ is given by 
$P(t)=\exp_{P_{0}}(t\dot{P}_{0})= \mu_{P_{0}}^{\exp}(t\dot{P}_{0})$. 
\end{remark}

Obviously, $\exp_{P}$ is smooth with
\begin{equation}
  \label{eq:109}
  \exp_{P}(0)=P
\end{equation}
and
\begin{equation}
  \label{eq:111}
  \DD\exp_{P}(0)=\id,
\end{equation}
as
\begin{equation}
  \label{eq:110}
  \begin{split}
   \DD\exp_{P}(0)\xi&=[[\xi,P],P]\\
&=\ad_{P}^{2}\xi\\&=\xi 
\qquad\qquad\text{for all}\; \xi\in
T_{P}\Gr_{m,n}. 
  \end{split}
\end{equation}
Such Riemannian normal coordinates
can be explicitly computed as follows. Given any $\Theta\in
\SO_{n}$ with
\begin{equation}
  \label{eq:67}
  P=\Theta^{\top}\begin{bmatrix}
      I_m & 0 \\
      0 &  0
    \end{bmatrix}\Theta
\end{equation}
we can write
\begin{equation}
  \label{eq:68}
  [\xi,P]=\Theta^{\top}\begin{bmatrix}
      0 & Z \\
      -Z^{\top} &  0
    \end{bmatrix}\Theta
\end{equation}
with $Z\in\R^{m\times(n-m)}$. Since
\begin{equation}
  \label{eq:69}
  [\xi,P]^{2m}=\Theta^{\top}\begin{bmatrix}
      (-ZZ^{\top})^{m}& 0 \\
      0 & (-Z^{\top}Z)^{m}
    \end{bmatrix}\Theta
\end{equation}
we obtain
\begin{equation}
  \label{eq:70}
  \e^{[\xi,P]}=\Theta^{\top}\e^{\left[\begin{smallmatrix}
      0 & Z \\
      -Z^{\top} &  0
    \end{smallmatrix}\right]}\Theta=
 \Theta^{\top} \begin{bmatrix}
    \cos\sqrt{ZZ^{\top}}&Z \ \frac{\sin\sqrt{Z^{\top}Z}}{\sqrt{Z^{\top}Z}}\\
-\frac{\sin\sqrt{Z^{\top}Z}}{\sqrt{Z^{\top}Z}} \ Z^{\top}&\cos\sqrt{Z^{\top}Z}
  \end{bmatrix}\Theta.
\end{equation}
Here, as usual, it is understood that
\begin{equation}
  \label{eq:72}
  \begin{split}
Z \ \frac{\sin\sqrt{Z^{\top}Z}}{\sqrt{Z^{\top}Z}}&=
Z(Z^{\top}Z)^{-\frac{1}{2}}
\sin(Z^{\top}Z)^{\frac{1}{2}}\\
&=Z\sum_{i=0}^{\infty}\frac{(-1)^{i}((Z^{\top}Z)^{\frac{1}{2}})^{2i}}{(2i+1)!}\\
&=Z\sum_{i=0}^{\infty}\frac{(-1)^{i}(Z^{\top}Z)^{i}}{(2i+1)!}\\
&=\phantom{Z}\sum_{i=0}^{\infty}\frac{(-1)^{i}(ZZ^{\top})^{i}}{(2i+1)!} \ Z\\
&=\frac{\sin\sqrt{ZZ^{\top}}}{\sqrt{ZZ^{\top}}} \ Z.
  \end{split}
\end{equation}
Therefore
\begin{equation}
  \begin{split}
  \label{eq:71}
  \exp_{P}(\xi)&=\Theta^{\top}\begin{bmatrix}
      \cos\sqrt{ZZ^{\top}} \\
       -\frac{\sin\sqrt{Z^{\top}Z}}{\sqrt{Z^{\top}Z}} \ Z^{\top}
\end{bmatrix}
      \begin{bmatrix}
      \cos\sqrt{ZZ^{\top}} &
-Z \ \frac{\sin\sqrt{Z^{\top}Z}}{\sqrt{Z^{\top}Z}} 
    \end{bmatrix}\Theta\\
    &=
    \Theta^{\top}\begin{bmatrix}
      \cos^2\sqrt{ZZ^{\top}} &-\cos\sqrt{ZZ^{\top}}\frac{\sin\sqrt{ZZ^{\top}}}{\sqrt{ZZ^{\top}}} \ Z\\
-Z^{\top}\frac{\sin\sqrt{ZZ^{\top}}}{\sqrt{ZZ^{\top}}}  \cos\sqrt{ZZ^{\top}} &
\sin^2\sqrt{Z^{\top}Z}
    \end{bmatrix}\Theta\\
    &=
    \Theta^{\top}\begin{bmatrix}
      \cos^2\sqrt{ZZ^{\top}} &-\sinc\left(2\sqrt{ZZ^{\top}}\right) \ Z\\
-Z^{\top}\sinc\left(2\sqrt{ZZ^{\top}}\right) &
\sin^2\sqrt{Z^{\top}Z}
    \end{bmatrix}\Theta\\
    &=
    \frac{1}{2}I_{n}+\Theta^{\top}\begin{bmatrix}
      \frac{1}{2}\cos\left(2\sqrt{ZZ^{\top}}\right) &-\sinc\left(2\sqrt{ZZ^{\top}}\right) \ Z\\
-Z^{\top}\sinc\left(2\sqrt{ZZ^{\top}}\right) &
-\frac{1}{2}\sin\left(2\sqrt{Z^{\top}Z}\right)
    \end{bmatrix}\Theta 
  \end{split}
\end{equation}

\subsubsection{QR-coordinates}
\label{sec:qr-coordinates}
We define $QR$-coordinates by the map
\begin{equation}\label{eq:QRco}
  \begin{split}
    \mu_{P}^{\mathrm{QR}}:T_{P}\Gr_{m,n} &\to \Gr_{m,n},\\
    \xi &\mapsto \left(I+[\xi,P]\right)_{\mathrm{Q}}P\left(\left(I+[\xi,P]\right)_{\mathrm{Q}}\right)^{\top}.
  \end{split}
  \end{equation}
  Here $M_{\mathrm{Q}}$ denotes the $Q$-factor in the $QR$-factorization
  $M=M_{\mathrm{Q}}M_{\mathrm{R}}$ of $M$.
Note that the matrix
\begin{equation}
  \label{eq:73}
  I + [\xi,P] = \Theta^{\top}\begin{bmatrix}
      I & Z \\
      -Z^{\top} &  I
    \end{bmatrix}\Theta
\end{equation}
is always invertible and therefore the $Q$-factor $(I+[\xi,P])_{\mathrm{Q}}\in O_{n}(\R)$
exists, and moreover, is unique if the diagonal entries of the upper triangular factor $R$ are chosen positive. 
From now on, we always choose the $R$-factor in this way. Actually, the determinant of the $Q$-factor,
\begin{equation}
  \label{eq:det}
  \det\left(I+[\xi,P]\right)_{\mathrm{Q}}=1,
\end{equation}
i.e., $\left(I+[\xi,P]\right)_{\mathrm{Q}}\in SO_{n}(\R)$, as it is easily checked that $\det
\left[\begin{smallmatrix}
  I&Z\\-Z^{\top}&I
\end{smallmatrix}\right]>0
$
always.
 
Moreover, by the smoothness of the
$QR$-factorization for general invertible matrices (follows from the
Gram-Schmidt procedure rather than from the usual algorithm via
Householder transformations), the map
$\mu_{P}^{\mathrm{QR}}$ is smooth on the tangent spaces
$T_{P}\Gr_{m,n}$ with
\begin{equation}
  \label{eq:112}
  \mu_{P}^{\mathrm{QR}}(0)=P.
\end{equation}
Now by \eqref{eq:112}
\begin{equation}
  \label{eq:113}
  \DD\mu_{P}^{\mathrm{QR}}(0):T_{P}\Gr_{m,n}\to T_{P}\Gr_{m,n}
\end{equation}
and a straightforward computation shows that
$\DD\mu_{P}^{\mathrm{QR}}(0)=\id$. In fact, by
differentiating the $QR$-factorization
\begin{equation}
  \label{eq:92}
  I+t[\xi,P]=Q(t)R(t),\ Q(0)=I,\ R(0)=I
\end{equation}
we obtain
\begin{equation}
  \label{eq:93}
  [\xi,P]=\dot{Q}R + Q\dot{R}
\end{equation}
and therefore at $t=0$
\begin{equation}
  \label{eq:94}
  [\xi,P]= \dot{Q}(0) + \dot{R}(0).
\end{equation}
But $[\xi,P]$ and $\dot{Q}(0)$ are skew-symmetric, while $\dot{R}(0)$
is upper triangular. Thus $\dot{R}(0)=0$ and therefore
\begin{equation}
  \label{eq:95}
  \left.\frac{\operatorname{d}}{\operatorname{d}t}(I+t[\xi,P])_{Q}\right|_{t=0}= \dot{Q}(0)=[\xi,P].
\end{equation}
This shows
\begin{equation}
  \label{eq:96}
  \DD\mu_{P}^{\mathrm{QR}}(0)\xi = [[\xi,P],P]=\xi
\end{equation}
for all $\xi \in T_{P}\Gr_{m,n}$, as claimed.

There exist explicit formulas for the $Q$ and $R$-factors in terms of Cholesky factors.
In fact, with
\begin{equation}\label{eq:P}
P=\Theta^{\top}\begin{bmatrix}
      I_m & 0 \\
      0 &  0
    \end{bmatrix}\Theta
\end{equation}
and since
\begin{equation}
  \label{eq:74}
  \Theta(I+[\xi,P])\Theta^{\top}=\begin{bmatrix}
      I_m & Z \\
      -Z^{\top} &  I_{n-m}
    \end{bmatrix}=:X
\end{equation}
satisfies
\begin{equation}
  \label{eq:75}
XX^{\top}=X^{\top}X=
\begin{bmatrix}
      I_m +ZZ^{\top}& 0 \\
      0 &  I_{n-m}+Z^{\top}Z
    \end{bmatrix}
\end{equation}
the $QR$-factorization of $X=
\left[\begin{smallmatrix}
   I_m & Z \\
      -Z^{\top} &  I_{n-m}
\end{smallmatrix}\right]=
X_{\mathrm{Q}}X_{\mathrm{R}}$ is obtained as
\begin{equation}
  \label{eq:XQ}
  X_{\mathrm{Q}}=
  \begin{bmatrix}
    R_{11}^{-1}&ZR_{22}^{-1}\\-Z^{\top}R_{11}^{-1}&R_{22}^{-1}
  \end{bmatrix}
\end{equation}
and 
\begin{equation}
  \label{eq:XR}
  X_{\mathrm{R}}=
  \begin{bmatrix}
    R_{11}&0\\0&R_{22}
  \end{bmatrix}.
\end{equation}
Here $R_{11}$ and $R_{22}$ are the unique Cholesky factors defined by 
\begin{equation}
  \begin{split}
    \label{Cholesky_R}
R_{11}^{\top}R_{11}&=I_{m}
+ZZ^{\top},\\
R_{22}^{\top}R_{22}&=I_{n-m}+Z^{\top}Z.
  \end{split}
\end{equation}
Note, that vanishing of the $12$-block of $X_{\mathrm{R}}$ follows from the invertibility of $R_{11}$ and $R_{22}$ and equation (\ref{eq:75}).



\subsubsection{Cayley coordinates}
\label{sec:cayley-coordinates}
Another possibility to introduce easily computable coordinates
utilizes the Cayley transform. For any 
skew-symmetric matrix $\Omega$ the {\em Cayley transform}
\begin{equation}
  \label{eq:23}
  \begin{split}
\Cay:\so_{n}&\to \SO_{n},\\
    \Omega &\to (2I+\Omega)(2I-\Omega)^{-1}
  \end{split}
\end{equation}
is smooth and satisfies $\DD\Cay(0)=\id$. The
Cayley coordinates are defined as
\begin{equation}\label{eq:QRco2}
  \begin{split}
    \mu_{P}^{\Cay}: T_{P}\Gr_{m,n} &\to \Gr_{m,n},\\
    \xi &\mapsto \Cay\left([\xi,P]\right)P\Cay\left(-[\xi,P]\right).
  \end{split}
  \end{equation}
The above mentioned property of the Cayley transform implies that
$\mu_{P}^{\Cay}$ is smooth and satisfies
\begin{equation}
  \label{eq:24}
  \begin{split}
    \mu_{P}^{\Cay}(0)&=P,\\
  \DD\mu_{P}^{\Cay}(0)\xi&=\DD\Cay(0)([\xi,P])P-P\DD\Cay(0)[\xi,P]=[[\xi,P],P]=\xi
\end{split}
\end{equation}
for all tangent vectors $\xi\in T_{P}\Gr_{m,n}$.
Moreover, $\mu_{P}^{\Cay}(\xi)$ is easily computed as follows. 
For
\begin{equation}
\label{eq:PZ}
P=\Theta^{\top}\begin{bmatrix}
      I_m & 0 \\
      0 &  0
    \end{bmatrix}\Theta \ , \quad
\xi=\Theta^{\top}\begin{bmatrix}
      0& -Z \\
      -Z^{\top}&  0
    \end{bmatrix}\Theta
\end{equation}
a straightforward computation shows, using Schur complements and properties of
the von Neumann series, that
\begin{equation*}
  \begin{split}
    \Cay\!\left(\!\begin{bmatrix}
      0& Z \\
      -Z^{\top}&  0
    \end{bmatrix}\!\right)
\!&=\!\begin{bmatrix}
      2I_{m}& Z \\
      -Z^{\top}&  2I_{n-m}
    \end{bmatrix}\!\!
\begin{bmatrix}
      2I_{m}& -Z \\
      Z^{\top}&  2I_{n-m}
    \end{bmatrix}^{-1}\\
\!&=\!\begin{bmatrix}
      2I_{m}& Z \\
      -Z^{\top}&  2I_{n-m}
    \end{bmatrix}\!\!
    \begin{bmatrix}
      \frac{1}{2}(I_m \!+\!\frac{1}{4}ZZ^{\top})^{-1}&
\frac{1}{4}Z(I_{n-m}\! +\!\frac{1}{4}Z^{\top}Z)^{-1}\\
-\frac{1}{4}Z^{\top}(I_{m}\! +\!\frac{1}{4}ZZ^{\top})^{-1}&
\frac{1}{2}(I_{n-m}\! +\!\frac{1}{4}Z^{\top}Z)^{-1}
    \end{bmatrix}\\
\!&=\!
    \begin{bmatrix}
      I_m -\frac{1}{4}ZZ^{\top}& Z \\
      -Z^{\top} &  I_{n-m}
      -\frac{1}{4}Z^{\top}Z
    \end{bmatrix}\!\!
    \begin{bmatrix}
      I_m \!+\!\frac{1}{4}ZZ^{\top}& 0 \\
      0&  I_{n-m}\! +\!\frac{1}{4}Z^{\top}Z
    \end{bmatrix}^{-1}.
  \end{split}
\end{equation*}
Therefore, 
\begin{equation} 
  \begin{split}
    \mu_{P}^{\Cay}(\xi)&=\Theta^{\top}
\Cay\left(\begin{bmatrix}
      0& Z \\
      -Z^{\top}&  0
    \end{bmatrix}\right)
    \begin{bmatrix}
      I_m & 0 \\
      0 &  0
    \end{bmatrix}
    \Cay\left(\begin{bmatrix}
      0& -Z \\
      Z^{\top}&  0
    \end{bmatrix}\right)\Theta\\
   & = \Theta^{\top} \begin{bmatrix}
      I_m -\frac{1}{4}ZZ^{\top}\\
      -Z^{\top}
    \end{bmatrix}
    \left(I_m +\frac{1}{4}ZZ^{\top}\right)^{-2}
    \begin{bmatrix}
      I_m -\frac{1}{4}ZZ^{\top}&-Z
    \end{bmatrix}\Theta.
  \end{split}
\end{equation}
Now
\begin{equation}
  \label{eq:81}
  \begin{bmatrix}
      I_m -\frac{1}{4}ZZ^{\top}\\
      -Z^{\top}
    \end{bmatrix}
    \left(I_m +\frac{1}{4}ZZ^{\top}\right)^{-1}
\end{equation}
is a basis matrix with orthonormal columns and therefore
$\mu_{P}^{\Cay}(\xi)$ is exactly the projection operator associated
with the linear subspace
\begin{equation}
  \label{eq:82}
  \Theta^{\top} \colspan
  \begin{bmatrix}
    I_m -\frac{1}{4}ZZ^{\top}\\
      -Z^{\top}
  \end{bmatrix}.
\end{equation}

\subsubsection{Approximation properties of parametrizations}
We have already shown that
\begin{equation}
  \label{eq:116}
  \DD\mu_{P}^{\exp}(0)=\DD\mu_{P}^{\Cay}(0)=
\DD\mu_{P}^{\mathrm{QR}}(0)=\id.
\end{equation}
Moreover, it holds
\begin{theorem}
  \label{sec:appr-prop-param-2ndAbleitung}
Let $P\in\Gr_{m,n}$ and $[\xi,P]$ 
as in \eqref{eq:67} and \eqref{eq:68}, respectively. Then
\begin{equation}
  \label{eq:120}
  \left.\textstyle{\frac{\operatorname{d}^{2}}{\operatorname{d}\varepsilon^{2}}}
\mu_{P}^{\exp}(\varepsilon\xi)\right|_{\varepsilon=0}\!\!\!\!\!=\!
\left.\textstyle{\frac{\operatorname{d}^{2}}{\operatorname{d}\varepsilon^{2}}}
\mu_{P}^{\Cay}(\varepsilon\xi)\right|_{\varepsilon=0}\!\!\!\!\!=\!
\left.\textstyle{\frac{\operatorname{d}^{2}}{\operatorname{d}\varepsilon^{2}}}
\mu_{P}^{\mathrm{QR}}(\varepsilon\xi)\right|_{\varepsilon=0}\!\!\!\!\!=\!
\Theta^{\top}\!
\begin{bmatrix}
  -2ZZ^{\top}&0\\0&2Z^{\top}\!Z
\end{bmatrix}
\!\Theta.
\end{equation}
Note, that the right hand side of \eqref{eq:120} is independent of the
choice of $\Theta$ in~\eqref{eq:67}.
\end{theorem}

\begin{proof}
  Taking derivatives yields
  \begin{equation}
    \label{eq:128}
    \begin{split}
      \left.\textstyle{\frac{\operatorname{d}^{2}}{\operatorname{d}\varepsilon^{2}}}
\mu_{P}^{\exp}(\varepsilon\xi)\right|_{\varepsilon=0}&= 
\left.\textstyle{\frac{\operatorname{d}^{2}}{\operatorname{d}\varepsilon^{2}}}
\e^{[\varepsilon\xi,P]}P\e^{-[\varepsilon\xi,P]}\right|_{\varepsilon=0}\\
&=[\xi,P]^{2}P+P[\xi,P]^{2}-
2[\xi,P]P[\xi,P]\\
&=\Theta^{\top}
\begin{bmatrix}
  -2ZZ^{\top}&0\\0&2Z^{\top}Z
\end{bmatrix}
\Theta
    \end{split}
  \end{equation}

From the theory of Pad{\'e} approximations it is well known that for all
matrices $X\in\so_{n}$ and $\varepsilon\in\R$
\begin{equation}
  \label{eq:129}
  \e^{\varepsilon X}=(2I_{n}+\varepsilon X)(2I_{n}-\varepsilon X)^{-1}+
\mathcal{O}(\varepsilon^{3}).
\end{equation}
Consequently,
\begin{equation}
  \label{eq:130}
  \left.\textstyle{\frac{\operatorname{d}^{2}}{\operatorname{d}\varepsilon^{2}}}
\mu_{P}^{\Cay}(\varepsilon\xi)\right|_{\varepsilon=0}=
\Theta^{\top}
\begin{bmatrix}
  -2ZZ^{\top}&0\\0&2Z^{\top}Z
\end{bmatrix}
\Theta
\end{equation}
holds as well. 

We now proceed with $\mu_{P}^{\mathrm{QR}}$. Let $\varepsilon\in\R$ be a parameter and let 
\begin{equation}
  \label{eq:Xe}
  X(\varepsilon):= \Theta(I+\varepsilon[\xi,P])\Theta^{\top}=
\begin{bmatrix}
      I_m & \varepsilon Z \\
      -\varepsilon Z^{\top} &  I_{n-m}
    \end{bmatrix}
\end{equation}
with $QR$-factorisation
\begin{equation}
  \label{eq:XeQR}
   X(\varepsilon)= X(\varepsilon)_{\mathrm{Q}} X(\varepsilon)_{\mathrm{R}},
\end{equation}
i.e.,
\begin{equation}
  \label{eq:XeQ}
  X(\varepsilon)_{\mathrm{Q}}= \begin{bmatrix}
    R_{11}^{-1}(\varepsilon)&\varepsilon ZR_{22}^{-1}(\varepsilon)\\-\varepsilon Z^{\top}R_{11}^{-1}(\varepsilon)&R_{22}^{-1}(\varepsilon)
  \end{bmatrix},
\end{equation}
where the Cholesky factors $R_{ii}$ are defined via
\begin{equation}
  \label{eq:Re}
  \begin{split}
    R_{11}^{\top}(\varepsilon)R_{11}(\varepsilon)&=I_{m}
+\varepsilon^{2}ZZ^{\top},\\
R_{22}^{\top}(\varepsilon)R_{22}(\varepsilon)&=I_{n-m}+\varepsilon^{2}Z^{\top}Z.
  \end{split}
\end{equation}
Obviously,
\begin{equation}
  \label{eq:Rii0}
  R_{11}(0)=I_{m}\qquad\text{and}\qquad R_{22}(0)=I_{n-m}.
\end{equation}
Therefore, taking the first order derivatives in (\ref{eq:Re}) and evaluating at $\varepsilon=0$ gives
\begin{equation}
  \label{eq:dotRii0}
  \begin{split}
    \dot{R}_{11}^{\top}(0)+\dot{R}_{11}(0)=0,\\
    \dot{R}_{22}^{\top}(0)+\dot{R}_{22}(0)=0,
  \end{split}
\end{equation}
which imply $\dot{R}_{11}(0)=0$ and $\dot{R}_{22}(0)=0$. Furthermore, taking second order derivatives at $\varepsilon=0$ and using (\ref{eq:Rii0}) and (\ref{eq:dotRii0}) gives
\begin{equation}
  \label{eq:ddotRii0}
  \begin{split}
    \ddot{R}_{11}^{\top}(0)+\ddot{R}_{11}(0)=2ZZ^{\top},\\
    \ddot{R}_{22}^{\top}(0)+\ddot{R}_{22}(0)=2Z^{\top}Z.
  \end{split}
\end{equation}
Using (\ref{eq:Rii0}), (\ref{eq:dotRii0}) and (\ref{eq:ddotRii0}) we compute the derivatives of the inverses as
\begin{equation}
  \label{eq:dotRinv}
  \begin{split}
    \left.\textstyle{\frac{\operatorname{d}}{\operatorname{d}\varepsilon}}R_{11}^{-1}(\varepsilon)\right|_{\varepsilon=0}&=0,\\
\left.\textstyle{\frac{\operatorname{d}}{\operatorname{d}\varepsilon}}R_{22}^{-1}(\varepsilon)\right|_{\varepsilon=0}&=0
  \end{split}
\end{equation}
and
\begin{equation}
  \label{eq:ddotRinv}
  \begin{split}
    \left.\textstyle{\frac{\operatorname{d}^{2}}{\operatorname{d}\varepsilon^{2}}}R_{11}^{-1}(\varepsilon)\right|_{\varepsilon=0}&=-\ddot{R}_{11}(0),\\
\left.\textstyle{\frac{\operatorname{d}^{2}}{\operatorname{d}\varepsilon^{2}}}R_{22}^{-1}(\varepsilon)\right|_{\varepsilon=0}&=-\ddot{R}_{22}(0).
  \end{split}
\end{equation}
Therefore, 
\begin{equation}
  \label{eq:XQmitdots}
  \begin{split}
    X(0)_{\mathrm{Q}}&=I,\\
\left.\textstyle{\frac{\operatorname{d}}{\operatorname{d}\varepsilon}}X(\varepsilon)_{\textrm{Q}}\right|_{\varepsilon=0}&=
\begin{bmatrix}
      0 &  Z \\
      - Z^{\top} &  0
    \end{bmatrix},\\
\left.\textstyle{\frac{\operatorname{d}^2}{\operatorname{d}\varepsilon^2}}X(\varepsilon)_{\textrm{Q}}\right|_{\varepsilon=0}&=
\begin{bmatrix}
       -\ddot{R}_{11}(0)& 0 \\
      0 &  -\ddot{R}_{22}(0)
    \end{bmatrix}
  \end{split}
\end{equation}
and finally,
\begin{equation}
  \begin{split}
   \left.\frac{\operatorname{d}^{2}}{\operatorname{d}\varepsilon^{2}}
\mu_{P}^{\mathrm{QR}}(\varepsilon\xi)\right|_{\varepsilon=0}\!\!\!&=
\Theta^{\top}\left(
\ddot{X}(0)_{\mathrm{Q}}
\left[\begin{smallmatrix}
  I_m&0\\0&0
\end{smallmatrix}\right]+
\left[\begin{smallmatrix}
  I_m&0\\0&0
\end{smallmatrix}\right]\ddot{X}^{\top}(0)_{\mathrm{Q}}+
2\dot{X}(0)_{\mathrm{Q}}
\left[\begin{smallmatrix}
  I_m&0\\0&0
\end{smallmatrix}\right]
\dot{X}^{\top}(0)_{\mathrm{Q}}
\right)\Theta\\
&=\Theta^{\top}
\begin{bmatrix}
  -\ddot{R}_{11}(0)-\ddot{R}^{\top}_{11}(0)&0\\0&2Z^{\top}Z
\end{bmatrix}\Theta\\
 &=\Theta^{\top}
\begin{bmatrix}
  -2ZZ^{\top}&0\\0&2Z^{\top}Z
\end{bmatrix}\Theta
 \end{split}
\end{equation}
as required.
\qed\end{proof}

\subsection{Gradients and Hessians}
Let $F:\Sym_{n}\to \R$ be a smooth function and
let $f:=F|_{\Gr_{m,n}}$ denote its restriction to the
Gra{\ss}mannian. Let $\nabla_{F}(P) \in \Sym_{n}$ be the gradient of $F$ in
$\Sym_{n}$ evaluated at $P$. Let $\operatorname{H}_{F}(P):
\Sym_{n}\times \Sym_{n}\to \R$ denote the
Hessian form of $F$ evaluated at $P$. We also consider 
$\Hess_{F}(P):\Sym_{n}\to \Sym_{n}$ as the corresponding linear map.
Gradient and Hessian are formed using the Euclidean (Frobenius) inner product on
$\Sym_{n}$. 

The next result computes the Riemannian gradient and Riemannian Hessian
of the restriction $f$ with respect to the induced Euclidean
Riemannian metric on the Gra{\ss}mannian (and thus also for the normal
Riemannian metric).

\begin{theorem}\label{the:Grass:gradients-hessians}
Let $f:\Gr_{m,n}\to \R$. 
The Riemannian gradient, $\mathfrak{grad}_{f}$, and the Riemannian Hessian operator
$\mathfrak{Hess}_{f}(P): T_{P}\Gr_{m,n}\to
T_{P}\Gr_{m,n}$ are given as
\begin{equation}
  \label{eq:40}
  \mathfrak{grad}_{f}(P) = \ad_{P}^{2}(\nabla_{F}(P)) = [P,[P,\nabla_{F}(P)]],
\end{equation}
and
\begin{equation}
  \label{eq:80}
  \mathfrak{Hess}_{f}(P)(\xi) = 
  \ad_{P}^{2}\Big(\operatorname{Hess}_{F}(P)(\xi)\Big) - 
  \ad_{P}\ad_{\nabla_{F}(P)}\xi,
\end{equation}
for all $\xi \in T_{P}\Gr_{m,n}$.
\end{theorem}

\begin{proof} The first part follows immediately from the well known 
fact, that the Riemannian gradient of $f$ coincides with the orthogonal
projection of $\nabla_{F}$ onto the tangent space
$T_{P}\Gr_{m,n}$. Since the orthogonal projection
operator onto the tangent space is given by $\ad_{P}^{2}$, this
proves the first claim. 

For the second part, consider a geodesic curve $P(t)$
with $\dot{P}(0)=\xi$. Then $\ddot{P}=-[\dot{P},[\dot{P},P]]$ and
therefore the Riemannian Hessian form is
\begin{equation}
    \label{eq:830}
  \begin{split}
 \mathfrak{H}_{f}(P(0))(\xi,\xi) &:= 
\left.\textstyle{\frac{\operatorname{d}^{2}(F\circ
    P)(t)}{\operatorname{d}t^2}}\right|_{t=0}\\
 &=\operatorname{H}_{F}(P(0))(\xi,\xi)+ \DD F(P(0)) \ \ddot{P}(0)\\
&=\operatorname{H}_{F}(P(0))(\xi,\xi)-\DD F(P(0)) \ [\xi,[\xi,P(0)]]\\
&= \operatorname{H}_{F}(P(0))(\xi,\xi)-\tr \Big(\nabla_{F}(P(0)) \ 
[\xi,[\xi,P(0)]]\Big).   
  \end{split}
\end{equation}
Thus by polarization
  \begin{equation}
    \label{eq:83}
    \begin{split}
   \mathfrak{H}_{f(P)}(\xi,\eta) &= \operatorname{H}_{F}(P)(\xi,\eta)-\half\tr
(\nabla_{F}(P) \ [\xi,[\eta,P]])-\frac{1}{2}\tr (\nabla_{F}(P) \ [\eta,[\xi,P]]))\\
&= \tr\left(\left(\Hess_{F}(P)(\xi) - \half[P,[\nabla_{F}(P),\xi]]
-\half[[\xi,P],\nabla_{F}(P)]\right)\eta \right)\\   
&= \tr\left(\left(\Hess_{F}(P)(\xi) - \half[P,[\nabla_{F}(P),\xi]]
-\half[\nabla_{F}(P),[P,\xi]\right)\eta \right)\\
&= \tr\Big(\ad_{P}^{2}\Big(\Hess_{F}(P)(\xi) -
\half[P,[\nabla_{F}(P),\xi]]
-\frac{1}{2}[\nabla_{F}(P),[P,\xi]]\Big)\eta \Big).   
    \end{split}
  \end{equation}
This implies that the Riemannian Hessian operator is given as
\begin{equation}
  \label{eq:84}
  \begin{split}
   \mathfrak{Hess}_{f}(P)(\xi)&= \ad_{P}^{2}\Big(\Hess_{F}(P)(\xi)-\half
\ad_{P}\ad_{\nabla_{F}(P)}\xi-\half\ad_{\nabla_{F}(P)}\ad_{P}\xi \Big)\\
&= \ad_{P}^{2}\Big(\Hess_{F}(P)(\xi)\Big) - \half
\ad_{P}^{3}\ad_{\nabla_{F}(P)}\xi-\half\ad_{P}^{2}\ad_{\nabla_{F}(P)}\ad_{P}\xi\\
&= \ad_{P}^{2}\Big(\Hess_{F}(P)(\xi)\Big) -  
\ad_{P}^{3}\ad_{\nabla_{F}(P)}\xi-\half\ad_{P}^{2}[\ad_{\nabla_{F}(P)},\ad_{P}]\xi\\ 
&= \ad_{P}^{2}\Big(\Hess_{F}(P)(\xi)\Big) -  \ad_{P}\ad_{\nabla_{F}(P)}\xi-
\half\ad_{P}^{2}[\ad_{\nabla_{F}(P)},\ad_{P}]\xi\\
&= \ad_{P}^{2}\Big(\Hess_{F}(P)(\xi)\Big) -  
\ad_{P}\ad_{\nabla_{F}(P)}\xi-\half\ad_{P}^{2}\ad_{[\nabla_{F}(P),P]}\xi. 
  \end{split}
\end{equation}
The result thus follows from the following lemma.
\qed\end{proof}
\begin{lemma}\label{lem:Grass-gradients-hessians}
For any tangent vector $\xi \in
T_{P}\Gr_{m,n}$ and any $A\in \Sym_{n}$ one has
\begin{equation}
  \label{eq:103}
  \ad_{P}\ad_{[A,P]}\xi=0.
\end{equation}
\end{lemma}

\begin{proof} Without loss of generality we can assume that
\begin{equation}
P=\begin{bmatrix}
      I_m & 0 \\
      0 &  0
    \end{bmatrix}, \quad
\xi=\begin{bmatrix}
      0& Z \\
      Z^{\top}&  0
    \end{bmatrix}, \quad
 A=\begin{bmatrix}
      A_{1}& A_{2} \\
      A_{2}^{\top}&  A_{3}
    \end{bmatrix}.
\end{equation}
Then
\begin{equation}
  \label{eq:104}
  [A,P]=\begin{bmatrix}
      0 & -A_{2} \\
      A_{2}^{\top}& 0
    \end{bmatrix}
\end{equation}
and
\begin{equation}
\ad_{[A,P]}\xi=\begin{bmatrix}
      -ZA_{2}^{\top}-A_{2}Z^{\top} & 0 \\
      0 &  A_{2}^{\top}Z+Z^{\top}A_2
    \end{bmatrix},
\end{equation}
from which $\ad_{P}\ad_{[A,P]}\xi=0$ follows by a straightforward
computation.
\qed\end{proof}

As a consequence we obtain the following formulas for the Riemannian gradient
and Riemannian Hessian operator of the Rayleigh quotient function.
\begin{corollary} \label{sec:gradients-hessians}Let $A\in \Sym_{n}$. The Riemannian gradient and Riemannian 
Hessian operator of the Rayleigh quotient function
\begin{equation}
  \label{eq:85}
  f: \Gr_{m,n}\to \R, \ f(P):= \tr(AP)
\end{equation}
are
  \begin{equation}
    \label{eq:86}
 \begin{split}   \mathfrak{grad}_{f}(P) &= [P,[P,A]],\\
    \mathfrak{Hess}_{f}(P)&= - \ad_{P}\circ \ad_{A},
 \end{split}
 \end{equation}
respectively.
\end{corollary}

As another example, let us consider the function
\begin{equation}
  \label{eq:91}
  \begin{split}
    F:\Sym_{n}&\to\R,\\
    P&\mapsto\|(I-P)AP\|^2=\tr(I-P)APA^{\top}
  \end{split}
\end{equation}
for an arbitrary, not necessarily symmetric matrix $A\in
\R^{n\times n}$. Note that the global minima of the restriction
$f:=F|_{\Gr_{m,n}}$ to the Gra{\ss}mannian are exactly the
projection operators corresponding to the $m$-dimensional
invariant subspaces of $A$. The gradient and Hessian operator on
$\Sym_{n}$ are computed as:
\begin{equation}
  \label{eq:149}
 \begin{split}      
\left.\frac{\operatorname{d}}{\operatorname{d}\varepsilon}
F(P+\varepsilon H)\right|_{\varepsilon=0}&=\tr(-HAPA^{\top}+(I-P)AHA^{\top})\\
&=\tr H(A^{\top}(I-P)A-APA^{\top}),\\
\left.\frac{\operatorname{d}^{2}}{\operatorname{d}\varepsilon^{2}}
F(P+\varepsilon H)\right|_{\varepsilon=0}&=-2\tr HAHA^{\top},
 \end{split}
\end{equation}
by polarisation for $H,K\in\Sym_{n}$
\begin{equation*}
\frac{1}{4}(-2\tr(H+K)A(H+K)A^{\top}\!\!\!+2(H-K)A(H-K)A^{\top})=
-\tr H(AKA^{\top}\!\!+A^{\top}KA)
\end{equation*}
consequently,
\begin{equation}
  \label{eq:150}
 \begin{split}      
 \nabla_{F}(P) &= A^{\top}(I-P)A - APA^{\top},\\
\Hess_{F}(P)(\xi) &= -A^{\top}\xi A - A\xi A^{\top}.   
 \end{split}
\end{equation}
This leads to the following explicit description of the Riemannian gradient
and Riemannian Hessian operators on the Gra{\ss}mannian.
\begin{corollary} Let $A\in \R^{n\times n}$. The Riemannian gradient and
  Riemannian 
Hessian operator of
  \begin{equation}
    \label{eq:88}
    f: \Gr_{m,n}\to \R, \ f(P):= \|(I-P)AP\|^2
  \end{equation}
are
\begin{equation}
  \label{eq:89}
  \begin{split}
 \mathfrak{grad}_{f}(P) &= [P,[P,A^{\top}A-A^{\top}PA-APA^{\top}]],\\
\mathfrak{Hess}_{f}(P)(\xi)&= -[P,[P,A^{\top}\xi A+A\xi
A^{\top}]]- [P,[A^{\top}A-A^{\top}PA-APA^{\top},\xi]],   
  \end{split}
\end{equation}
respectively.
\end{corollary}

\section{Geometry of the Lagrange Gra{\ss}mannian}
In this section we develop an analogous theory for the manifold of
Lagrangian subspaces in $\R^{2n}$. Thus we consider the
{\em Lagrange Gra{\ss}mann manifold}  $\LG_{n}$, consisting
of all $n$-dimensional Lagrangian subspaces of $\R^{2n}$
with respect to the standard symplectic form
\begin{equation}
  \label{eq:105}
  J:=\begin{bmatrix}
      0 & I_{n} \\
      -I_{n} &  0
    \end{bmatrix}.
\end{equation}
Recall, that an $n$-dimensional subspace $V\subset \R^{2n}$ is called Lagrangian, if
\begin{equation}
  \label{eq:106}
  v^{\top}Jv=0
\end{equation}
for all $v\in V$.
Instead of interpreting the elements of the Lagrange Gra{\ss}mann 
manifold as maximal isotropic subspaces, we prefer to view them
in an equivalent way as a certain subclass of symmetric projection
operators. Note that, if $P$ is the symmetric projection operator
onto an $n$-dimensional linear subspace $V$, then the condition
$PJP=0$ is equivalent to $V$ being Lagrangian. Thus we define the
{\em Lagrange Gra{\ss}mannian}
\begin{equation}
  \label{eq:25}
  \LG_{n}:= \{P \in \Sym_{2n}\ | \ P^2=P,\tr P=n, PJP=0\}
\end{equation}
as the manifold of rank $n$ symmetric projection operators of
$\R^{2n}$, satisfying the Lagrangian subspace condition $PJP=0$.
Note, that $\LG_{n}$ is a compact, connected submanifold of
the Gra{\ss}mannian $\Gr_{n,2n}$. In order to obtain a
deeper understanding of the geometry of this set, we observe that
$\LG_{n}$ is a homogeneous space for the action of the
orthogonal symplectic group. Let
\begin{equation}
  \label{eq:14}
  \GL_{n}:=\{X\in\R^{n\times n}|\det X\ne 0\}
\end{equation}
and let 
\begin{equation}
  \label{eq:26}
  \OSp_{2n}:= \{T\in \mathrm{GL}_{2n}\ | \
T^{\top}J T = J,\ T \in \SO_{2n} \}
\end{equation}
denote the Lie
group of orthogonal symplectic transformations. Let
\begin{equation}
  \label{eq:27}
  \osp_{2n}=\{X\in \so_{2n}\ | \
X^{\top}J + JX =0\}
\end{equation}
denote the
associated Lie algebra of skew-symmetric Hamiltonian $2n\times
2n$-matrices. Thus the elements of $\osp_{2n}$ are
exactly the real $2n\times 2n$-matrices $T$ of the form
\begin{equation}
  \label{eq:28}
  T:=\begin{bmatrix}
      A & -B \\
      B &  A
    \end{bmatrix}
\end{equation}
defined by the condition, that $A+\imath B\in \mathfrak{u}_{n}$,
i.e. $A+\imath B$ is 
skew-Hermitian, i.e., $A\in\so_{n}$ and $B\in \Sym_{n}$, where 
\begin{equation}
  \label{eq:29}
  \mathfrak{u}_{n}:=\{X\in\C^{n\times n}|X^{\ast}=-X\}
\end{equation}
and the asterisk symbol denotes complex conjugate transpose and 
$\imath:=\sqrt{-1}$.
Similarly, the elements of $\OSp_{2n}$
are the real $2n\times 2n$-matrices $\xi$ of the form
\begin{equation}
  \label{eq:30}
  \xi:=\e^{\left[\begin{smallmatrix}
      X & -Y \\
      Y &  X
    \end{smallmatrix}\right]}\quad\text{satisfying}\quad X\in\so_{n},\;Y\in\Sym_{n}.
\end{equation}
In particular,
$\OSp_{2n}$ is isomorphic to the unitary
group
\begin{equation}
  \label{eq:31}
  \U_{n}:=\{X\in\C^{n\times n}|X^{\ast}X=I_{n}\}.
\end{equation}
The orthogonal symplectic group 
$\OSp_{2n}$ acts transitively on $\LG_{n}$ via
\begin{equation}\label{eq:1}
  \begin{split}
    \sigma:\OSp_{2n}\times\LG_{n}&\to
    \LG_{n},\\
    (T,P)&\mapsto T^{\top}PT,
  \end{split}
\end{equation}
with the stabilizer subgroup of
\begin{equation}
  \label{eq:32}
  P:=\begin{bmatrix}
      I & 0 \\
      0 &  0
    \end{bmatrix}\in \LG_{n}
\end{equation}
given as the set of all block-diagonal matrices
\begin{equation}
  \label{eq:33}
  T=\begin{bmatrix}
      A & 0 \\
      0 & A
    \end{bmatrix}\quad\text{with}\quad A\in \mathrm{O}_{n}. 
\end{equation}
Therefore $\LG_{n}$ is a
homogeneous space that can be identified with
$\mathrm{U}_{n}/\mathrm{O}_{n}$.

\begin{theorem}\label{thm:basic}
\begin{enumerate}
\item[(a)]The Lagrange Gra{\ss}mannian $\LG_{n}$ is a
smooth, compact submanifold of $\Sym_{2n}$ of dimension
$\frac{n(n+1)}{2}$. \item[(b)] The tangent space of
$\LG_{n}$ at an element $P\in \LG_{n}$ is given as
\begin{equation}
  \label{eq:34}
  T_P\LG_{n}=\{[P,\Omega]\ | \ \Omega\in \mathfrak{osp}_{2n}\}.
\end{equation}
\end{enumerate}
\end{theorem}

Since the tangent space
$T_{P}\LG_{n}\subset \Sym_{2n}$ is a subset
of $\Sym_{2n}$, we can define the normal space at $P$ to be
the vector space
\begin{equation}
  \label{eq:35}
 N_{P}\LG_{n}= (T_{P}\LG_{n})^{\perp} := \{X\in \Sym_{2n}\ | \ \tr(XY)=0
\ \text{for all}\ Y\in T_{P}\LG_{n}\}.   
\end{equation}
\begin{proposition} 
\label{sec:geom-lagr-grassm}
Let $P\in \LG_{n}$ be arbitrary.
\begin{enumerate}
\item 
The normal subspace in $\Sym_{2n}$ at $P$ is given as
\begin{equation}
  \label{eq:36}
  N_{P}\LG_{n}=\left\{X-\textstyle{\frac{1}{2}}\ad_{P}^{2}(JXJ+X)\ | 
\ X\in \Sym_{2n}\right\}.
\end{equation}
\item
The linear map
\begin{equation}
  \label{eq:45}
  \begin{split}
    \pi: \Sym_{2n}&\to \Sym_{2n},\\
 X&\mapsto \textstyle{\frac{1}{2}}[P,[P,JXJ+X]]
  \end{split}
\end{equation}
is a self-adjoint projection operator onto
$T_{P}\LG_{n}$ with kernel $N_{P}\LG_{n}$.
\end{enumerate}
\end{proposition}

\begin{proof} 
To prove the first statement let $\Omega\in\osp_{2n}$, $X\in\Sym_{2n}$ and 
$P\in\LG_{n}$ be arbitrary. Then
for $[P,\Omega]\in T_{P}\LG_{n}$ and 
$X-\frac{1}{2}\ad_{P}^{2}(JXJ+X)\in N_{p}\LG_{n}$ we
have
\begin{equation}
  \label{eq:37}
  \begin{split}
    \tr\left([P,\Omega]\left(X-\textstyle{\frac{1}{2}}\ad_{P}^2(JXJ+X)\right)\right)&=
\tr\left(\Omega\ \left(\textstyle{\frac{1}{2}}\ad_{P}^3(JXJ+X)-\ad_{p}X\right)\right)\\
 &=\tr\left(\Omega\left(\textstyle{\frac{1}{2}}\ad_{P}(JXJ+X)-\ad_{P}X)\right)\right)\\
 &=\textstyle{\frac{1}{2}}\tr\left(\Omega\ \ad_{P}(JXJ-X)\right)\\
 &=\tr\left(\Omega P(JXJ-X)\right)\\
 &=\tr\left(X\Omega(JPJ-P)\right)\\
&= 0,
  \end{split}
\end{equation}
where we have used Lemma~\ref{lem:adP3}, $\Omega$ being skew-symmetric and
Hamiltonian, and the easily verified identity 
\begin{equation}
  \label{eq:15}
  JPJ-P=-I_{2n},\;\text{for all}\; P\in\LG_{n}.
\end{equation}
By \eqref{eq:37}, 
$T_{P}\LG_{n}$ and $N_{P}\LG_{n}$ are orthogonal subspaces of $\Sym_{2n}$
with respect to the Frobenius inner product. Analogously to 
Proposition~\ref{prop:normalGR} we now see that 
$\Sym_{2n}=T_{P}\LG_{n}\oplus N_{P}\LG_{n}$ holds true as well:

Every $P\in\LG_{n}$ can be written as
\begin{equation}
  \label{eq:38}
  P=Q^{\top}
  \begin{bmatrix}
    I_{n}&0\\0&0
  \end{bmatrix}Q
\end{equation}
for some $Q\in\OSp_{2n}$. Note that for all $\Omega\in\osp_{2n}$, and for all 
$X,S\in\Sym_{2n}$
\begin{equation}
  \label{eq:39}
\begin{split}    
\tr\left(S[P,\Omega]\right)&=\tr\left(QSQ^{\top}
\left[
  \left[\begin{smallmatrix}
  I_{n}&0\\0&0
  \end{smallmatrix}\right],Q\Omega Q^{\top}
\right]\right),\\
  \tr\!\left(S\!\left(\!X\!\!-\!\!\textstyle{\frac{1}{2}}\!
\ad_{P}^2(JXJ\!+\!\!X)\!\right)\!\right)&=
\tr\!\big(QSQ^{\top}\!\!\big(QXQ^{\top}\!\!\!-
\!\textstyle{\frac{1}{2}}\!\ad^{2}_{\left[\begin{smallmatrix}
      I_{n} & 0 \\ 0 & 0
  \end{smallmatrix}\right]}(JQXQ^{\top}\!\!\!J\!+\!QXQ^{\top})\big)\big).
\end{split}
\end{equation}
By \eqref{eq:39} and $\Sym_{2n}\to Q(\Sym_{2n})Q^{\top}$ 
being an isomorphism, 
without loss of generality we might assume that
\begin{equation}
  \label{eq:41}
  P=\begin{bmatrix}
  I_{n}&0\\0&0
  \end{bmatrix}.
\end{equation}
Assume there exists an $S\in\Sym_{2n}$ being orthogonal to both
subspaces. We will show the implication
\begin{equation}
  \label{eq:42}
  \left.
    \begin{aligned}
      \tr\left(S\ad_{P}\Omega\right)&=0\quad\text{for all}\quad\Omega\in\osp_{2n}\\
    \tr\left(S(X-\textstyle{\frac{1}{2}}\ad_{P}^{2}(JXJ+X)\right)&=
0\quad\text{for all}\quad
    X\in\Sym_{2n}
\end{aligned}\right\}\Longrightarrow\;S=0.
\end{equation}
Partition $S\in\Sym_{2n}$ as
\begin{equation}
  \label{eq:43}
  S=
  \begin{bmatrix}
    S_{11}&S_{12}\\S_{12}^{\top}&S_{22}
  \end{bmatrix}
\end{equation}
then
\begin{equation}
  \label{eq:44}
  \tr\left(S\ad_{P}\Omega\right)=0\quad\text{for
    all}\quad\Omega\in\osp_{2n}\;
\Longleftrightarrow\;S_{12}\in\so_{n},
\end{equation}
where we have used the symmetry of the $(12)-$block of $\Omega$, see
\eqref{eq:27} and \eqref{eq:28}. Moreover, 
\begin{equation}
  \label{eq:46}
  \begin{split}
     \tr\left(S(X-\half\ad_{P}^{2}(JXJ+X)\right)&=
\tr\left(SX-\half(\ad_{P}^{2}S)(JXJ+X)\right)\\
&=\tr\left(
  \begin{bmatrix}
  S_{11}&\frac{S_{12}-S_{12}^{\top}}{2}\\\frac{S_{12}^{\top}-S_{12}}{2}&S_{22}
  \end{bmatrix}X
\right).
  \end{split}
\end{equation}
Therefore,
\begin{equation}
  \label{eq:47}
  \tr\left(S(X-\half\ad_{P}^{2}(JXJ+X)\right)\!=\!0\;\text{for all}\;
    X\in\Sym_{2n}\;\Longleftrightarrow\; 
\left\{
  \begin{aligned}
    &S_{12}\in\Sym_{n},\\
&S_{11}=S_{22}=0.
  \end{aligned}\right.
\end{equation}
Together with \eqref{eq:44} we conclude $S=0$, i.e., 
$\Sym_{2n}=T_{P}\LG_{n}\oplus N_{P}\LG_{n}$ as required.

Now we prove the second claim. By the same reasoning as above we again 
might assume that 
\begin{equation}
  \label{eq:410}
  P=\begin{bmatrix}
  I_{n}&0\\0&0
  \end{bmatrix}.
\end{equation}
Let
\begin{equation}
  \label{eq:51}
  X=\begin{bmatrix}
    X_{11}&X_{12}\\X_{12}^{\top}&X_{22}
  \end{bmatrix}.
\end{equation}
Since
\begin{equation}
  \label{eq:52}
  \pi(X)=\half\ad_{P}^{2}(JXJ+X)=\begin{bmatrix}
   0&\frac{X_{12}+X_{12}^{\top}}{2}\\\frac{X_{12}+X_{12}^{\top}}{2}&0
  \end{bmatrix},
\end{equation}
we see that $\pi^{2}(X)=\pi(X)$
and moreover $\im \pi = T_{P}\LG_{n}.$ For any 
$X-\half\ad_{P}^{2}(JXJ+X) \in N_{P}\LG_{n}$ 
we have
\begin{equation}
  \label{eq:180}
  \pi\left(X-\half\ad_{P}^{2}(JXJ+X)\right)=\pi(X)-\pi^{2}(X)=0,
\end{equation}
and by 
counting dimensions $\ker \pi = N_{P}\LG_{n}$. Finally, for
all $X,Y\in\Sym_{2n}$ and by using~\eqref{eq:52} 
\begin{equation}
\label{eq:53}
\begin{split}
\langle \pi(X),Y\rangle &=\tr\left(\half\ad_{P}^{2}(JXJ+X)Y\right)\\
&=\half\tr\left(\begin{bmatrix}
   0&X_{12}+X_{12}^{\top}\\X_{12}+X_{12}^{\top}&0
  \end{bmatrix}\begin{bmatrix}
   Y_{11}&Y_{12}\\Y_{12}^{\top}&Y_{22}
  \end{bmatrix}\right)\\
&=\half\tr\left(\begin{bmatrix}
   X_{11}&X_{12}\\X_{12}^{\top}&X_{22}
  \end{bmatrix}\begin{bmatrix}
   0&Y_{12}+Y_{12}^{\top}\\Y_{12}+Y_{12}^{\top}&0
  \end{bmatrix}\right)\\
&=\tr\left(X\ \half\ad_{P}^{2}(JYJ+Y)\right)\\
&=\langle X,\pi(Y)\rangle.  
\end{split}
\end{equation}
Thus $\pi$ is self-adjoint and the result follows.
\qed\end{proof}

Fortunately, the discussion of Riemannian metrics carries directly over from
the Gra{\ss}mannian case to the case of the Lagrange Gra{\ss}mannian. We
therefore omit the proof.

Consider the surjective linear map
\begin{equation}
  \label{eq:190}
  \begin{split}
    \ad_P: \osp_{2n} &\to T_{P}\LG_{n},\\ \Omega&\mapsto [P,\Omega]
  \end{split}
\end{equation}
with kernel
\begin{equation}
  \label{eq:1020}
  \ker\ad_P = \{\Omega \in \osp_{2n} \ | \ P\Omega = \Omega P\}.
\end{equation}
We regard $\osp_{2n}$ as an inner product space, endowed
with the Frobenius inner product $\langle \Omega_1, \Omega_2
\rangle = \tr(\Omega^{\top}\Omega_2)$. Then $\ad_P$ induces an
isomorphism of vector spaces
\begin{equation}
  \label{eq:550}
  \widehat{\ad}_P: (\ker \ad_P)^{\perp} \to T_{P}\LG_{n}
\end{equation}
and therefore induces an isometry of inner product spaces, by
defining an inner product on $T_{P}\LG_{n}$ via
\begin{equation}
  \label{eq:560}
  \langle\langle X,Y \rangle\rangle _{P} := -\tr(\widehat{\ad}_{P}^{-1}(X)
\widehat{\ad}_P^{-1}(Y))
\end{equation}
called the normal Riemannian metric. 

\begin{proposition}
The Euclidean and normal Riemannian metrics on
the Lagrange Gra{\ss}mannian $\LG_{n}$ coincide, i.e. for all
$P\in \LG_{n}$ and for all $X,Y \in
T_{P}\LG_{n}$ we have
\begin{equation}
  \label{eq:570}
  \tr(X^{\top}Y)=-\tr\left(\widehat{\ad}_{P}^{-1}(X)\ \widehat{\ad}_P^{-1}(Y)\right).
\end{equation}
\end{proposition}

Since a solution to~\eqref{eq:geodesic1} with an initial value
$P_{0}\in\LG_{n}$ and $\dot P_{0}\in T_{P_0}\LG_{n}$ is fully
contained in $\LG_{n}$, and since geodesics are unique, the geodesics of 
$\LG_{n}$ are also described by that equation.

\begin{theorem} 
\label{theo:GeoLagrGrass}
The geodesics of $\LG_{n}$ are
exactly the solutions of the second order differential equation
\begin{equation}
\label{eq:geodesic10}
\ddot{P} + [\dot{P},[\dot{P},P]] = 0.
\end{equation}
The unique geodesic $P(t)$ with initial conditions $P(0)=P_0 \in
\LG_{n}$, $\dot{P}(0)=\dot{P}_0 \in
T_{P_0}\LG_{n}$ is given by
\begin{equation}
\label{eq:geodesic20}
P(t) = \e^{t[\dot{P}_0,P_0]}P_0\e^{-t[\dot{P}_0,P_0]}.
\end{equation}
\end{theorem}

We now consider {\em local parametrizations} for the
Lagrange Gra{\ss}mannian as well
\begin{equation}
  \label{eq:660}
 \mu_P : T_{P}\LG_{n}\to \LG_{n} 
\end{equation}
satisfying
\begin{equation}
  \label{eq:1080}
  \mu_{P}(0)=P\;\text{and} \;\DD\mu_{P}(0)=\id.
\end{equation}

\subsection{Parametrizations and coordinates for the Lagrange
  Gra{\ss}\-man\-nian}
\label{sec:LG-coord}
\subsubsection{Riemannian normal coordinates}
\label{sec:LG-riem-norm-coord}
As before Riemannian normal coordinates are defined through
\begin{equation}\label{eq:LGRMcoords}
  \begin{split}
    \exp_{P}: T_{P}\LG_{n} &\to \LG_{n},\\
    \exp_{P}(\xi) &= \e^{[\xi,P]}P\e^{-[\xi,P]}.
  \end{split}
\end{equation}
Given any $\Theta\in
\OSp_{2n}$ with
\begin{equation}
  \label{eq:670}
  P=\Theta^{\top}\begin{bmatrix}
      I_n & 0 \\
      0 &  0
    \end{bmatrix}\Theta
\end{equation}
and
\begin{equation}
  \label{eq:680}
  [\xi,P]=\Theta^{\top}\begin{bmatrix}
      0 & Z \\
      -Z &  0
    \end{bmatrix}\Theta,\quad Z\in\Sym_{n}.
\end{equation}
We obtain
\begin{equation}
  \label{eq:710}
  \exp_{P}(\xi)=\Theta^{\top}\begin{bmatrix}
      \phantom{-}\cos Z \\
       -\sin Z
\end{bmatrix}
      \begin{bmatrix}
      \cos Z &
-\sin Z 
    \end{bmatrix}\Theta.
\end{equation}

\subsubsection{QR-coordinates}
\label{sec:LGqr-coordinates}
Let $P$ and $[\xi,P]$ as in \eqref{eq:670} and \eqref{eq:680}. 
We define smooth $QR$-coordinates by the map
\begin{equation}\label{eq:LGQRmodco}
  \begin{split}
    \mu_{P}^{\mathrm{QR}}: T_{P}\LG_{n} &\to \LG_{n},\\
    \xi &\mapsto
    \left(I+[\xi,P]\right)_{\mathrm{Q}}
    P
    \left(I+[\xi,P]\right)_{\mathrm{Q}}^{\top}.
  \end{split}
\end{equation}
Analogous to Section~\ref{sec:qr-coordinates} we define
\begin{equation}
  \label{eq:7700}
  \begin{bmatrix}
      I_n & Z \\
      -Z &  I_{n}
    \end{bmatrix}_{\mathrm{Q}}:= 
\begin{bmatrix}
      R^{-1}&ZR^{-1}\\-ZR^{-1}&R^{-1}
    \end{bmatrix}
\end{equation}
where $R$ denotes the unique Cholesky factor that solves $R^{\top}R=I+Z^2$.
Note that the $\mathrm{Q}$-factor in \eqref{eq:7700} is orthogonal and
symplectic.
The above map~\eqref{eq:LGQRmodco} is therefore well defined.

\subsubsection{Cayley coordinates}
\label{sec:LGcayley-coordinates}
Let $P$ and $[\xi,P]$ be as in \eqref{eq:670} and \eqref{eq:680}. For any
skew-symmetric Hamiltonian matrix $\Omega$ the {\em Cayley transform}
\begin{equation}
  \label{eq:230}
  \begin{split}
\Cay:\osp_{2n}&\to \OSp_{2n},\\
    \Omega &\to (2I+\Omega)(2I-\Omega)^{-1}
  \end{split}
\end{equation}
is smooth and satisfies $\DD\Cay(0)=\id$. The
Cayley coordinates are defined as
\begin{equation}\label{eq:LGCco2}
  \begin{split}
    \mu_{P}^{\Cay}: T_{P}\LG_{n} &\to \LG_{n},\\
    \xi &\mapsto \Cay\left([\xi,P]\right)P\Cay\left(-[\xi,P]\right).
  \end{split}
  \end{equation}
Therefore,
\begin{equation} 
  \begin{split}
    \mu_{P}^{\Cay}(\xi)
   & = \Theta^{\top} \begin{bmatrix}
      I_n -\frac{1}{4}Z^{2}\\
      -Z
    \end{bmatrix}
    \left(I_n +\frac{1}{4}Z^{2}\right)^{-2}
    \begin{bmatrix}
      I_n -\frac{1}{4}Z^{2}&-Z
    \end{bmatrix}\Theta.
  \end{split}
\end{equation}

\subsection{Gradients and Hessians}
Let $F:\Sym_{n}\to \R$ be a smooth function and
let $f:=F|_{\LG_{n}}$ denote its restriction to the
Lagrange Gra{\ss}mannian. Let $\nabla_{F}(P) \in \Sym_{n}$ be the gradient of $F$ in
$\Sym_{n}$ evaluated at $P$. Let $\operatorname{H}_{F}(P):
\Sym_{n}\times \Sym_{n}\to \R$ denote the
Hessian form of $F$ evaluated at $P$. We also consider 
$\Hess_{F}(P):\Sym_{n}\to \Sym_{n}$ as the corresponding linear map.
Gradient and Hessian are formed using the Euclidean (Frobenius) inner product on
$\Sym_{n}$. 

The next result computes the Riemannian gradient and Riemannian Hessian
of the restriction $f$ with respect to the induced Euclidean
Riemannian metric on the Lagrange Gra{\ss}mannian (and thus also for the normal
Riemannian metric).

\begin{theorem}
\label{the:LG-gradients-hessians}
Let $f:\LG_{n}\to \R$ and consider the orthogonal 
projection operator $\pi$ as defined by \eqref{eq:45}. 
The Riemannian gradient, $\mathfrak{grad}_{f}$, and the Riemannian Hessian operator
$\mathfrak{Hess}_{f}(P): T_{P}\LG_{n}\to
T_{P}\LG_{n}$ are given as
\begin{equation}
  \label{eq:400}
  \mathfrak{grad}_{f}(P) = \pi(\nabla_{F}(P))=
  \half\ad_{P}^{2}\Big(J(\nabla_{F}(P))J+\nabla_{F}(P)\Big),
\end{equation}
and
\begin{equation}
  \begin{split}    
    \label{eq:800}
    \mathfrak{Hess}_{f}(P)(\xi)& = 
    \half\ad_{P}^{2}\Big(J\Big(\operatorname{Hess}_{F}(P)(\xi)\Big)J+
    \operatorname{Hess}_{F}(P)(\xi)\Big)\\
    &\quad - \half\ad_{P}^{2}\Big(J\Big(\ad_{P}
    \ad_{\nabla_{F}(P)}\xi\Big)J+\ad_{P}\ad_{\nabla_{F}(P)}\xi\Big),
  \end{split}
\end{equation}
for all $\xi \in T_{P}\Gr_{m,n}$.
\end{theorem}

\begin{proof} The first part follows again from the
fact, that the Riemannian gradient of $f$ coincides with the orthogonal
projection of $\nabla_{F}$ onto the tangent space
$T_{P}\LG_{n}$. 

Analogously to the proof of Theorem~\ref{the:Grass:gradients-hessians} 
\begin{equation}
    \label{eq:8300}
 \mathfrak{H}_{f}(P(0))(\xi,\xi) = 
\operatorname{H}_{F}(P(0))(\xi,\xi)-\tr \Big(\nabla_{F}(P(0)) \ 
[\xi,[\xi,P(0)]]\Big).   
\end{equation}
Thus by polarization
  \begin{equation}
    \label{eq:8311}
    \begin{split}
   \mathfrak{H}_{f}(P)(\xi,\eta) 
&= \tr\left(\left(\Hess_{F}(P)(\xi) -
    \half[P,[\nabla_{F}(P),\xi]]-\half[\nabla_{F}(P),[P,\xi]\right)\eta
\right)\\ 
&= \tr\Big(\pi\Big(\Hess_{F}(P)(\xi) -
\half[P,[\nabla_{F}(P),\xi]]-\frac{1}{2}[\nabla_{F}(P),[P,\xi]\Big)\eta
\Big).    
    \end{split}
  \end{equation}
This implies that the Riemannian Hessian operator is given as
\begin{equation}
  \label{eq:8410}
  \begin{split}
   \mathfrak{Hess}_{f}(P)(\xi)&= \pi\Big(\Hess_{F}(P)(\xi)-\half
\ad_{P}\ad_{\nabla_{F}(P)}\xi-\half\ad_{\nabla_{F}(P)}\ad_{P}\xi \Big)\\
&= \pi\Big(\Hess_{F}(P)(\xi)
-\ad_{P}\ad_{\nabla_{F}(P)}\xi
-\half\ad_{[\nabla_{F}(P),P]}\xi \Big)\\
&= \pi\Big(\Hess_{F}(P)(\xi)\Big)
-\pi\Big(\ad_{P}\ad_{\nabla_{F}(P)}\xi\Big)
-\pi\Big(\half\ad_{[\nabla_{F}(P),P]}\xi \Big).
  \end{split}
\end{equation}
Together with Lemma~\ref{lem:Grass-gradients-hessians} the
Lemma~\ref{lem:Lagrange-gradients-hessians} below implies
\begin{equation}
  \label{eq:114}
  \pi\Big(\half\ad_{[\nabla_{F}(P),P]}\xi \Big)=0.
\end{equation}
The result follows.
\qed\end{proof}
\begin{lemma}\label{lem:Lagrange-gradients-hessians} For any tangent vector $\xi \in
T_{P}\LG_{n}$ and any $A\in \Sym_{n}$ one has
\begin{equation}
  \label{eq:1030}
  \ad_{P}(J(\ad_{[A,P]}\xi)J)=0.
\end{equation}
\end{lemma}

\begin{proof} Without loss of generality we can assume that
\begin{equation}
P=\begin{bmatrix}
      I_m & 0 \\
      0 &  0
    \end{bmatrix}, \quad
\xi=\begin{bmatrix}
      0& Z \\
      Z&  0
    \end{bmatrix}, \quad
Z=Z^{\top},\quad
 A=\begin{bmatrix}
      A_{1}& A_{2} \\
      A_{2}^{\top}&  A_{3}
    \end{bmatrix}.
\end{equation}
Then
\begin{equation}
J(\ad_{[A,P]}\xi)J=\begin{bmatrix}
     -A_{2}^{\top}Z -ZA_{2} & 0 \\
      0 &  ZA_{2}^{\top}+A_2 Z
    \end{bmatrix},
\end{equation}
from which $\ad_{P}(J(\ad_{[A,P]}\xi)J)=0$ follows by a straightforward
computation.
\qed\end{proof}

As a consequence we obtain the following formulas for the Riemannian gradient
and Riemannian Hessian operator of the Rayleigh quotient function used
to compute
the $n-$dimensional dominant eigenspace of a real symmetric Hamiltonian
$(2n\times 2n)-$matrix.

Consider the set of real symmetric Hamiltonian
$(2n\times 2n)-$matrices $\mathfrak{p}_{2n}$
\begin{equation}
  \label{eq:122}
  \mathfrak{p}_{2n}:=\left\{H\in\Sym_{2n}\ \left| \ H=
    \begin{bmatrix}
     S&T\\T&-S   
    \end{bmatrix}\right.,\;T,S\in\Sym_{n}.
\right\}
\end{equation}
Note that 
\begin{equation}
  \label{eq:123}
  JKJ=K\qquad\text{for all}\quad K\in\mathfrak{p}_{2n}.
\end{equation}
Moreover, from the theory of Cartan decompositions, see e.g. \cite{knapp:96},  
the following commutator relations are well known
\begin{equation}
  \label{eq:124}
  [\mathfrak{p}_{2n},\mathfrak{p}_{2n}]\subset\osp_{2n},\quad 
[\mathfrak{p}_{2n},\osp_{2n}]\subset\mathfrak{p}_{2n}
\end{equation}
together with the isomorphisms
\begin{equation}
  \label{eq:125}
  Q^{\top}(\osp_{2n})Q\cong\osp_{2n}, \quad
  Q^{\top}(\mathfrak{p}_{2n})Q\cong\mathfrak{p}_{2n}
\quad\text{for all}\quad Q\in\OSp_{2n}.
\end{equation}
\begin{corollary} Given $H\in\mathfrak{p}_{2n}$
Let 
\begin{equation}
\label{eq:115}
  F:\Sym_{2n}\to \R,\quad P\mapsto\tr(HP),
\end{equation}
with restriction
\begin{equation}
  \label{eq:117}
  f:=F|_{\LG_{n}}.
\end{equation}
The Riemannian gradient and Riemannian 
Hessian operator are
  \begin{equation}
    \label{eq:860}
 \begin{split}   \mathfrak{grad}_{f}(P) &= [P,[P,A]],\\
    \mathfrak{Hess}_{f}(P)&= - \ad_{P}\circ \ad_{A},
 \end{split}
 \end{equation}
respectively.
\end{corollary}
\begin{proof}
Because the function $F$ is linear the Euclidean gradient is simply
  \begin{equation}
    \label{eq:118}
    \nabla_{F}(P)=H,
  \end{equation}
and the Euclidean Hessian operator vanishes
\begin{equation}
  \label{eq:119}
  \Hess_{F}(P)=0.
\end{equation}
We therefore get for the Riemannian gradient using~\eqref{eq:123} and  
Theorem~\ref{the:LG-gradients-hessians}
\begin{equation}
  \label{eq:121}
  \begin{split}
    \mathfrak{grad}_{f}(P)&=\half\ad_{P}^{2}\Big(J(\nabla_{F}(P))J+\nabla_{F}(P)\Big)\\
&=\half\ad_{P}^{2}\Big(JHJ+H\Big)\\
&=[P,[P,H]].
  \end{split}
\end{equation}
For the Riemannian Hessian operator we need some preparation. Let
$\xi=[P,\Omega]\in T_{P}\LG_{n}$ be arbitrary, i.e., $\Omega\in\osp_{2n}$ is
arbitrary. Then there exists a $Q\in\OSp_{2n}$ such that 
\begin{equation}
  \label{eq:126}
  \xi=Q^{\top}\left[
    \begin{bmatrix}
    I_{n}&0\\0&0
    \end{bmatrix},Q\Omega Q^{\top}
\right]Q.
\end{equation}
But the commutator in~\eqref{eq:126} is an element of $\mathfrak{p}_{2n}$ and
therefore by~\eqref{eq:125} the same holds true for $\xi$ independent of
$P$. Consequently, $\ad_{H}\xi\in\osp_{2n}$ and 
$\ad_{P}\ad_{H}\xi\in\mathfrak{p}_{2n}$, and finally
$J(\ad_{P}\ad_{H}\xi)J=\ad_{P}\ad_{H}\xi$. The formula for the Riemannian
Hessian operator is now easily verified
\begin{equation}
  \label{eq:127}
  \begin{split}
    \mathfrak{Hess}_{f}(P)(\xi)&=
-\half\ad_{P}^{2}\Big(J(\ad_{P}\ad_{\nabla_{F}(P)}\xi)J
+\ad_{P}\ad_{\nabla_{F}(P)}\xi\Big)\\ 
&=-\half\ad_{P}^{2}\Big(J(\ad_{P}\ad_{H}\xi)J+\ad_{P}\ad_{H}\xi\Big)\\
&=-\ad_{P}\ad_{H}\xi.
  \end{split}
\end{equation}
\qed\end{proof}

\section{Newton's method}
In the following we propose a class of Newton-like algorithms to compute
a nondegenerate critical point of a smooth cost function
$f:\Gr_{m,n}\to\R$. Local quadratic convergence of the proposed algorithm will
be established.
Parts of this section are based on the conference paper \cite{huep:05a}. 

\subsection{The Euclidean case}
\label{sec:newtons-method-rn}
Let $f:\R^{n}\to\R$ be a smooth function and let $x^{*}\in\R^{n}$
be a nondegenerate critical point of $f$, i.e. the Hessian operator 
$\Hess_{f}(x^{*})$ is
invertible. Newton's method for $f$ is the iteration
\begin{equation}
  \label{eq:1z}
  x_{0}\in\R^{n}, \ x_{k+1}=N_{f}(x_{k}):=x_{k}-\left(\Hess_{f}(x_{k})
\right)^{-1}\nabla_{f}(x_{k}).
\end{equation}
Note that the iteration~\eqref{eq:1z} is only defined if $\Hess_{f}(x_{k})$ is
invertible for all $k\in\N_{0}=\N\cup\{0\}$. However, since $f$ is smooth, there exists an
open neighborhood of $x^{*}$ in which the Hessian operator is invertible.

It is well known that the point sequence $\{x_{k}\}_{k\in\N_{0}}$ generated
by~\eqref{eq:1z} is defined and converges locally quadratically to $x^{*}$
provided that $x_{0}$ is sufficiently close to $x^{*}$. For more information
see e.g. \cite{lue:84}. 

\subsection{The Gra{\ss}mannian case}
Let $\{\mu_{P}\}_{P\in\Gr_{m,n}}$ be a family of local
parametri\-zations of $\Gr_{m,n}$.
Let $P^{*}\in\Gr_{m,n}$ be a nondegenerate critical
point of the smooth function $f:\Gr_{m,n}\to\R$. 
If there exists an open neighborhood
$U\subset\Gr_{m,n}$ of $P^{*}$ and a smooth map
\begin{equation*}
  \mu:U\times\R^{m(n-m)}\longrightarrow\Gr_{m,n}
\end{equation*}
such that $\mu(P,x)=\mu_{P}(x)$ for all $P\in U$ and $x\in\R^{m(n-m)}$ we will call 
$\{\mu_{P}\}_{P\in\Gr_{m,n}}$ a \emph{locally smooth family of
  parametrizations around $P^{*}$}. 
Let $\{\mu_{P}\}_{P\in
  \Gr_{m,n}}$ and $\{\nu_{P}\}_{P\in \Gr_{m,n}}$ be two locally smooth families of
parametrizations around $P^{*}$.
Consider the following iteration on $\Gr_{m,n}$
\begin{equation}
  \label{eq:2z}
  P_{0}\in \Gr_{m,n}, \ P_{k+1}=\nu_{P_{k}}\left(N_{f\circ\mu_{P_{k}}}(0)\right)
\end{equation}
where $N_{f\circ\mu_{P}}$ is defined in~\eqref{eq:1z}. The following theorem
is an adaptation from \cite{huep:05a}, where it is stated and proved
for arbitrary smooth manifolds.

\begin{theorem}\label{th:main-result}
  Under the condition
  \begin{equation}
    \label{eq:3z}
    \DD\mu_{P^{*}}(0)=\DD\nu_{P^{*}}(0)
  \end{equation}
  there exists an open neighborhood $V\subset \Gr_{m,n}$ of $P^{*}$ such that the
  point sequence $\{P_{k}\}_{k\in\N_{0}}$ generated by~\eqref{eq:2z} converges
  quadratically to $P^{*}$ provided $P_{0}\in V$.
\end{theorem}

\begin{proof}
Let $\mu,\nu:U\times\R^{m(n-m)}\to\Gr_{m,n}$ be smooth and such that
$\mu(P,x)=\mu_{P}(x)$ and $\nu(P,x)=\nu_{P}(x)$ for all $P\in U$ and
$x\in\R^{m(n-m)}$, where $U$ is a neighborhood of $P^{*}$.

The derivative of the algorithm map
\begin{equation}
  \begin{split}
    \label{eq:135}
  s:\Gr_{m,n}&\to \Gr_{m,n},\\ 
P&\mapsto \nu\left(P,-\left(\Hess_{f\circ\mu}(P,0)\right)^{-1}
\nabla_{f\circ\mu}(P,0)\right)
  \end{split}
\end{equation}
at $P^{*}$ is the linear map 
\begin{equation}
   \label{eq:139}
  \DD s(P^{*}):T_{P^{*}}\Gr_{m,n}\to T_{P^{*}}\Gr_{m,n}
\end{equation}
defined by
\begin{equation}
  \label{eq:140}
  \begin{split}
  & \DD s(P^{*})h= 
\DD_{1}\nu\left(P^{*},
-\left(\Hess_{f\circ\mu}(P^{*},0)\right)^{-1}
\nabla_{f\circ\mu}(P^{*},0)\right)h\\
&+\DD_{2}\nu\left(P^{*},
-\left(\Hess_{f\circ\mu}(P^{*},0)\right)^{-1}
\nabla_{f\circ\mu}(P^{*},0)\right)\\
&\quad\cdot
\left(-\DD_{P}\left(\left(\Hess_{f\circ\mu}(P^{*},0)\right)^{-1}
\nabla_{f\circ\mu}(P^{*},0)\right)h\right).
  \end{split}
\end{equation}
Here $\DD_{i}(\cdot)h$ denotes the derivative of $(\cdot)$ with respect to the
$i-$th argument in direction $h$, whereas, by abuse of notation, $\DD_{P}$
denotes the differential operator to compute the derivative with respect to 
the argument $P$.

The first summand on the right side in \eqref{eq:140} is easily computed as
\begin{equation}
\label{eq:13z}
  \begin{split}
\DD_{1}\nu\left(P^{*},
-\left(\Hess_{f\circ\mu}(P^{*},0)\right)^{-1}
\nabla_{f\circ\mu}(P^{*},0)\right)h&=\DD_{1}\nu(P^{*},0)h=h,
  \end{split}
\end{equation}
which is true because $P^{*}$ is a critical point of $f$ and the
gradient therefore vanishes. 

The second summand in \eqref{eq:140} consists of two terms due to the chain
rule. We first compute the left term giving 
\begin{equation}
  \label{eq:141}
  \begin{split}
    \DD_{2}\nu\left(P^{*},
-\left(\Hess_{f\circ\mu}(P^{*},0)\right)^{-1}
\nabla_{f\circ\mu}(P^{*},0)\right)= \DD_{2}\nu(P^{*},0)
  \end{split}
\end{equation}
because $P^{*}$ is critical. The evaluation of the right term is more
involved.
\begin{multline}
  \label{eq:142}
    -\DD_{P}\left(\left(\Hess_{f\circ\mu}(P^{*},0)\right)^{-1}
\nabla_{f\circ\mu}(P^{*},0)\right)h=\\
\shoveleft{\hspace*{1cm}
 -\left(\DD_{P}\left(\left(\Hess_{f\circ\mu}(P^{*},0)\right)^{-1}\right)h\right)
 \underbrace{\nabla_{f\circ\mu}(P^{*},0)}_{=0}}\\
\shoveright{-\left(\Hess_{f\circ\mu}(P^{*},0)\right)^{-1}\DD_{P}
 \left(\nabla_{f\circ\mu}(P^{*},0)\right)h\hspace*{1cm}}\\
=-\left(\Hess_{f\circ\mu}(P^{*},0)\right)^{-1}\DD_{P}
\left(\nabla_{f\circ\mu}(P^{*},0)\right)h.
\end{multline}
By the definition of the gradient one has for any $x\in\R^{m(n-m)}$
\begin{equation}
  \label{eq:143}
  \langle\nabla_{f\circ\mu}(P,0),x\rangle=
\DD f(P)\cdot\DD_{2}\mu(P,0)\cdot x
\end{equation}
and therefore using the critical point condition and the definition of the
Hessian operator in terms of second derivatives
\begin{equation}
  \label{eq:144}
  \begin{split}
      \left\langle\DD_{P}\left(\nabla_{f\circ\mu}(P^{*},0)\right)h,x\right\rangle&=
\DD^{2} f(P^{*})\cdot\left(\DD_{2}\mu(P^{*},0)\cdot x,h\right)\\
&\quad+\underbrace{\DD f(P^{*})}_{=0}\DD_{1}\left(\DD_{2}\mu(P^{*},0)\cdot
      x\right)h\\
&=\DD^{2} f(P^{*})\cdot\left(\DD_{2}\mu(P^{*},0)\cdot x,h\right)\\
&=\DD^{2} f(P^{*})\cdot\left(\DD_{2}\mu(P^{*},0)
\cdot x,\DD_{2}\mu(P^{*},0)\left(\DD_{2}\mu(P^{*},0)\right)^{-1}h\right)\\
&=\left\langle\Hess_{f\circ\mu}(P^{*},0)(x),
\left(\DD_{2}\mu(P^{*},0)\right)^{-1}h\right\rangle\\
&=\left\langle\Hess_{f\circ\mu}(P^{*},0)
\left(\DD_{2}\mu(P^{*},0)\right)^{-1}h,x\right\rangle
  \end{split}
\end{equation}
We now can conclude
\begin{equation}
  \label{eq:145}
  \DD_{P}\left(\nabla_{f\circ\mu}(P^{*},0)\right)h= \Hess_{f\circ\mu}(P^{*},0)
\left(\DD_{2}\mu(P^{*},0)\right)^{-1}h
\end{equation}
which in turn implies that
\begin{equation}
  \label{eq:146}
  \left(\Hess_{f\circ\mu}(P^{*},0)\right)^{-1} \DD_{P}\left(\nabla_{f\circ\mu}(P^{*},0)\right)h= 
\left(\DD_{2}\mu(P^{*},0)\right)^{-1}h.
\end{equation}
Summarizing our computations we have shown that
\begin{equation}
  \begin{split}
  \label{eq:147}
   \DD
   s(P^{*})h&=h-\DD_{2}\nu(P^{*},0)\left(\DD_{2}\mu(P^{*},0)\right)^{-1}h\\
&=0.
  \end{split}
\end{equation}
Consider now a local representation of $s$ in coordinate charts around
$P^{*}$ and $s(P^{*})=P^{*}$. Let $\|\cdot\|$ denote any norm in the
local coordinate space. By abuse of notation we will still speak of
$s$, $P^{*}$ and so on in reference to their local coordinate representations. 
Using a Taylor expansion of $s$ around $P^{*}$, there exists a neighborhood
$\overline{V}_{P^{*}}$ of $P^{*}$ such that the estimate
\begin{equation}
  \label{eq:148}
  \|s(P)-P^{*}\|\leq
\sup_{Q\in\overline{V}_{P^{*}}}\|\DD^{2}s(Q)\|\cdot\|P-P^{*}\|^{2}
\end{equation}
holds for all $P\in \overline{V}_{P^{*}}$.
Therefore, the subset $U\subset\overline{V}_{P^{*}}$
\begin{equation*}
U:=\{P\in\overline{V}_{P^{*}}\;|\;\sup_{Q\in\overline{V}_{P^{*}}}\|\DD^{2}s(Q)\|\cdot\|P-P^{*}\|<1
\}
\end{equation*}
is a neighborhood of $P^{*}$ that is
invariant under $s$, and hence remains invariant under the iterations of $s$. This
completes the proof of local quadratic convergence of the algorithm. 
\qed\end{proof}

A few remarks are in order.
Geometrically, the iteration~\eqref{eq:2z} does the following. The current
iteration point $P_{k}$ is pulled back to Euclidean space via the local
parametrization $\mu_{P_{k}}$ around $P_{k}$. Then one Euclidean Newton step
is performed for the function expressed in local coordinates, followed by a
projection back onto the Gra{\ss}mannian using the local parametrization $\nu_{P_{k}}$
around $P_{k}$.

For the special choice 
$\{\mu_{p}\}_{p\in  M}=\{\nu_{p}\}_{p\in M}$ both Riemannian normal coordinates
(cf. Section~\ref{sec:riem-norm-coord}), our iteration~\eqref{eq:2z} is precisely the so-called
Newton method along geodesics of D. Gabay \cite{gabay:82}, more recently also referred
to as the intrinsic Newton method. This follows from the lemma below.
\begin{lemma}
  \label{sec:grassmannian-case}
Let $f:\Gr_{m,n}\to \R$ be a smooth function.
For all $P\in\Gr_{m,n}$, $\xi\in T_{P}\Gr_{m,n}$ and any
$\mathtt{type}\in\{\exp,\mathrm{QR},\Cay\}$ we have
\begin{equation}
  \label{eq:50}
  \nabla_{f\circ\mu_{P}^{\mathtt{type}}}(0)=
  \mathfrak{grad}_{f}(P)
\end{equation}
and
\begin{equation}
\label{eq:152}
\operatorname{Hess}_{f\circ\mu_{P}^{\mathtt{type}}}(0)(\xi)=
\mathfrak{Hess}_{f}(P)(\xi).
\end{equation}
\end{lemma}

\begin{proof}
  By Remark~\ref{sec:riem-norm-coord-2} the unique geodesic through
  a point $P\in\Gr_{m,n}$ in direction $\xi\in T_{P}\Gr_{m,n}$ is
  given by $P(t)=\mu_{P}^{\exp}(t\xi)$ and hence
  \begin{equation}
    \label{eq:151}
    \langle\langle \mathfrak{grad}_{f}(P),\xi \rangle\rangle_{P}=
    \left.\frac{\operatorname{d}}{\operatorname{d}\varepsilon}(f\circ\mu_{P}^{\exp})
    (\varepsilon)\right|_{\varepsilon=0}=
    \langle\nabla_{f\circ\mu_{P}^{\exp}}(0),\xi\rangle,
  \end{equation}
  which by Proposition~\ref{prop:equalmetrics} implies \eqref{eq:50}
  for $\mathtt{type}=\exp$. The result for
  $\mathtt{type}=\mathrm{QR}$ and $\mathtt{type}=\Cay$ then
  follows from~\eqref{eq:116}.
  By the same line of arguments~\eqref{eq:152} follows from
  \begin{equation}
    \label{eq:153}
    \langle\langle\mathfrak{Hess}_{f}(P)(\xi),\xi\rangle\rangle_{P}=
    \left.\frac{\operatorname{d}^{2}}{\operatorname{d}\varepsilon^{2}}(f\circ\mu_{P}^{\exp})
    (\varepsilon)\right|_{\varepsilon=0}=
    \langle\operatorname{Hess}_{f\circ\mu_{P}^{\mathtt{\exp}}}(0)(\xi),\xi\rangle
  \end{equation}
  and Theorem~\ref{sec:appr-prop-param-2ndAbleitung}.
\qed\end{proof}

\subsection{The Lagrange Gra{\ss}mannian case}
All the above results carry over literally to Newton-like algorithms
on the Lagrange Gra{\ss}mannian by substituting the respective formulas.

\subsection{Algorithms}
\label{sec:algorithm}
We conclude by presenting several specific instances of the resulting
algorithms.
We discuss the case of smooth functions $F:\Sym_{n}\to \R$ with
restriction $f:=F|_{\Gr_{m,n}}$ to the Gra{\ss}mannian, as well as the
special cases of the Rayleigh quotient 
function on the Gra{\ss}mannian and the Lagrange Gra{\ss}mannian.
Furthermore, we consider the previously introduced nonlinear trace function for invariant
subspace computations on the Gra{\ss}mannian.
In all cases we choose $\{\mu_{P}\}$ as the Riemannian normal coordinates and
$\{\nu_{P}\}$ as the QR-coordinates, 
see Sections~\ref{sec:param-coord-grassm} and~\ref{sec:LG-coord}.

Recall that our convergence result requires the Hessian of the
restricted function to be nondegenerate at the critical point.

We first formulate a preliminary form of the algorithm we are
interested in.

\hspace{-5mm}\fbox{\begin{minipage}[t]{125mm}
\noindent Step 1.
\begin{itemize}
\item[] Pick a rank $m$ symmetric projection operator of $\R^{n}$, 
$P_{0}\in\Gr_{m,n}$, and
set $j=0$.
\end{itemize}
Step 2.
\begin{itemize}
\item[] Solve
  \begin{equation*}
    \ad_{P_{j}}^{2}\operatorname{Hess}_{F}(P_{j})(\ad_{P_{j}}\Omega_{j})
    -\ad_{P_{j}}\ad_{\nabla_{F}(P_{j})}\ad_{P_{j}}\Omega_{j}=
    -\ad_{P_{j}}^{2}\nabla_{F}(P_{j})
  \end{equation*}
for $\Omega_{j}\in\so(n)$.
\end{itemize}
Step 3.
\begin{itemize}
\item[] Solve
  \begin{equation*}
    P_{j}=\Theta_{j}^{\top}\begin{bmatrix}
      I_m & 0 \\
      0 &  0
    \end{bmatrix}\Theta_{j}
  \end{equation*}
  for $\Theta_{j}\in\SO_{n}$.
\end{itemize}
Step 4.
\begin{itemize}
\item[] Compute
  \begin{equation*}
    P_{j+1}=\Theta_{j}^{\top}
    \left(\Theta_{j}(I-\ad_{P_{j}}^{2}\Omega_{j})\Theta_{j}^{\top}\right)_{\mathrm{Q}}
    \!\!\!\!\Theta_{j} P_{j}\Theta_{j}^{\top}
    \left(\Theta_{j}(I-\ad_{P_{j}}^{2}\Omega_{j})\Theta_{j}^{\top}\right)_{\mathrm{Q}}^{\top}
    \!\!\!\!\Theta_{j}.
  \end{equation*}
\end{itemize}
Step 5.
\begin{itemize}
\item[] Set $j=j+1$ and goto Step 2.
\end{itemize}
\end{minipage}}\\

Here the expressions in Step 2 result from applying
Lemma~\ref{sec:grassmannian-case} and 
Theorem~\ref{the:Grass:gradients-hessians} and using the
representation $\xi=[P,\Omega]$, $\Omega\in\so(n)$ for an element 
$\xi\in T_{P}\Gr_{m,n}$.

Inspecting the above algorithm, it is evident that it can be
rewritten as an iteration in the $\Theta_{j}\in\SO_{n}$ as follows.

\hspace{-5mm}\fbox{\begin{minipage}[t]{125mm}
\noindent Step 1.
\begin{itemize}
\item[] Pick an orthogonal matrix $\Theta_{0}\in\SO_{n}$ corresponding
  to 
  \begin{equation*}
    P_{0}=\Theta_{0}^{\top}\begin{bmatrix}
      I_m & 0 \\
      0 &  0
    \end{bmatrix}\Theta_{0}\in\Gr_{m,n},
  \end{equation*}
  and set $j=0$.
\end{itemize}
Step 2.
\begin{itemize}
\item[] Solve
  \begin{equation*}
    \ad_{P_{j}}^{2}\operatorname{Hess}_{F}(P_{j})(\ad_{P_{j}}\Omega_{j})
    -\ad_{P_{j}}\ad_{\nabla_{F}(P_{j})}\ad_{P_{j}}\Omega_{j}=
    -\ad_{P_{j}}^{2}\nabla_{F}(P_{j})
  \end{equation*}
for $\Omega_{j}\in\so(n)$.
\end{itemize}
Step 3.
\begin{itemize}
\item[] Compute
  \begin{equation*}
    \Theta_{j+1}^{\top}=\Theta_{j}^{\top}
    \left(\Theta_{j}(I-\ad_{P_{j}}^{2}\Omega_{j})\Theta_{j}^{\top}\right)_{\mathrm{Q}}
    \text{and}\quad
    P_{j+1}=\Theta_{j+1}^{\top}\begin{bmatrix}
      I_m & 0 \\
      0 &  0
    \end{bmatrix}\Theta_{j+1}
  \end{equation*}
\end{itemize}
Step 4.
\begin{itemize}
\item[] Set $j=j+1$ and goto Step 2.
\end{itemize}
\end{minipage}}\\

Note that by equation~\eqref{eq:74} the term
\begin{equation}
  \label{eq:154}
  \Theta_{j}(I-\ad_{P_{j}}^{2}\Omega_{j})\Theta_{j}^{\top}=\begin{bmatrix}
      I_m & Z_{j} \\
      -Z_{j}^{\top} &  I_{n-m}
    \end{bmatrix}
\end{equation}
of which we have to compute a QR-factorization in Step 3 has a nice
block structure that can be exploited to get an efficient implementation.
Locally quadratic convergence is guaranteed as long as the specific
QR-factorization used is differentiable,
cf. Section~\ref{sec:qr-coordinates}. 

For specific functions the necessary computations might drastically
simplify, as the following two examples show for the Rayleigh quotient
function $F:\Sym_{n}\to \R,\ F(P)=\tr(AP)$, $A\in\Sym_{n}$, 
cf. Corollary~\ref{sec:gradients-hessians}.

\subsubsection{Rayleigh quotient on the Gra{\ss}mannian}
\label{alg:rayleigh-grassmann}

The equation we have to solve for $\Omega_{j}\in\so(n)$ in Step 2
becomes
\begin{equation}
  -\ad_{P_{j}}\ad_{A}\ad_{P_{j}}\Omega_{j}=
  -\ad_{P_{j}}^{2}A
\end{equation}
which is equivalent to
\begin{equation}
  \Theta_{j}(\ad_{P_{j}}\ad_{A}\ad_{P_{j}}\Omega_{j})\Theta_{j}^{\top}=
  \Theta_{j}(\ad_{P_{j}}^{2}A)\Theta_{j}^{\top}
\end{equation}
and, using 
\begin{equation}
  \label{eq:155}
  P_{j}=\Theta_{j}^{\top}\begin{bmatrix}
    I_m & 0 \\
    0 &  0
  \end{bmatrix}\Theta_{j},
\end{equation}
is equivalent to solving
\begin{equation}
  \label{eq:156}
  \ad_{\left[\begin{smallmatrix}I_{m} & 0 \\ 0 & 0\end{smallmatrix}\right]}
  \ad_{\Theta_{j}A\Theta_{j}^{\top}}
  \ad_{\left[\begin{smallmatrix}I_{m} & 0 \\ 0 & 0\end{smallmatrix}\right]}
  \begin{bmatrix}0 & Z_{j} \\ -Z_{j}^{\top} & 0\end{bmatrix}=
  \ad_{\left[\begin{smallmatrix}I_{m} & 0 \\ 0 & 0\end{smallmatrix}\right]}^{2}
  (\Theta_{j}A\Theta_{j}^{\top})
\end{equation}
for $Z_{j}\in\R^{m\times(n-m)}$. Denoting
\begin{equation}
  \label{eq:157}
  \Theta_{j}A\Theta_{j}^{\top}=
  \begin{bmatrix}A_{11} & A_{12} \\ A_{12}^{\top} & A_{22}\end{bmatrix}
\end{equation}
we actually have to solve the Sylvester equation
\begin{equation}
  \label{eq:158}
  A_{11}Z_{j}-Z_{j}A_{22}=A_{12}.
\end{equation}
The resulting algorithm is exactly the algorithm presented
in~\cite{huep:05a}. 

\hspace{-5mm}\fbox{\begin{minipage}[t]{125mm}
    \begin{center}
      \textbf{Algorithm 1: Rayleigh quotient on the Gra{\ss}mannian.}
    \end{center}
\noindent Step 1.
\begin{itemize}
\item[] Pick an orthogonal matrix $\Theta_{0}\in\SO_{n}$ corresponding
  to 
  \begin{equation*}
    P_{0}=\Theta_{0}^{\top}\begin{bmatrix}
      I_m & 0 \\
      0 &  0
    \end{bmatrix}\Theta_{0}\in\Gr_{m,n},
  \end{equation*}
  and set $j=0$.
\end{itemize}
Step 2.
\begin{itemize}
\item[] Compute
  \begin{equation*}
    \begin{bmatrix}A_{11} & A_{12} \\ A_{12}^{\top} & A_{22}\end{bmatrix}=
  \Theta_{j}A\Theta_{j}^{\top}.
\end{equation*}
\end{itemize}
Step 3.
\begin{itemize}
\item[] Solve the Sylvester equation
  \begin{equation*}
    A_{11}Z_{j}-Z_{j}A_{22}=A_{12}.
  \end{equation*}
for $Z_{j}\in\R^{m\times(n-m)}$.
\end{itemize}
Step 4.
\begin{itemize}
\item[] Compute
  \begin{equation*}
    \Theta_{j+1}^{\top}=\Theta_{j}^{\top}
    \begin{bmatrix}
      I_m & Z_{j} \\
      -Z_{j}^{\top} &  I_{n-m}
    \end{bmatrix}_{\mathrm{Q}}
    \quad\text{and}\quad
    P_{j+1}=\Theta_{j+1}^{\top}\begin{bmatrix}
      I_m & 0 \\
      0 &  0
    \end{bmatrix}\Theta_{j+1}
  \end{equation*}
\end{itemize}
Step 4.
\begin{itemize}
\item[] Set $j=j+1$ and goto Step 2.
\end{itemize}
\end{minipage}}\\

Since the global maximum of the Rayleigh quotient function on the
Gra{\ss}mannian is a nondegenerate critical point, provided that there
is a spectral gap after the $m$th largest eigenvalue of $A$, we immediately get
the following result.

\begin{theorem}
  For almost all matrices $A\in\Sym_{n}$
  Algorithm 1 converges locally quadratically to the projector onto the
  $m$-dimensional dominant eigenspace.
\end{theorem}

\subsubsection{Rayleigh quotient on the Lagrange Gra{\ss}mannian}
\label{alg:rayleigh-lagrange}
The corresponding problem of optimizing the Rayleigh quotient function
over the Lagrange Gra{\ss}mann manifold can be treated completely analogous to
the approach above. Thus let $A$ denote a real symmetric Hamiltonian matrix of
size $2n \times 2n$. The Newton algorithm for optimizing the trace function
$\tr (AP)$ over $\LG_{n}$ then is as follows. 

\hspace{-5mm}\fbox{\begin{minipage}[t]{125mm}
    \begin{center}
      \textbf{Algorithm 2: Rayleigh quotient on the Lagrange Gra{\ss}mannian.}
    \end{center}
\noindent Step 1.
\begin{itemize}
\item[] Pick an orthogonal matrix $\Theta_{0}\in\SO_{2n}$ corresponding
  to 
  \begin{equation*}
    P_{0}=\Theta_{0}^{\top}\begin{bmatrix}
      I_n & 0 \\
      0 &  0
    \end{bmatrix}\Theta_{0}\in\LG_{n},
  \end{equation*}
  and set $j=0$.
\end{itemize}
Step 2.
\begin{itemize}
\item[] Compute
  \begin{equation*}
    \begin{bmatrix}A_{11} & A_{12} \\ A_{12} & -A_{11}\end{bmatrix}=
  \Theta_{j}A\Theta_{j}^{\top}.
\end{equation*}
\end{itemize}
Step 3.
\begin{itemize}
\item[] Solve the Lyapunov equation
  \begin{equation*}
    A_{11}Z_{j}+Z_{j}A_{11}=A_{12}.
  \end{equation*}
for the symmetric matrix $Z_{j}\in\Sym_{n}$.
\end{itemize}
Step 4.
\begin{itemize}
\item[] Compute
  \begin{equation*}
    \Theta_{j+1}^{\top}=\Theta_{j}^{\top}
    \begin{bmatrix}
      I_n & Z_{j} \\
      -Z_{j} &  I_{n}
    \end{bmatrix}_{\mathrm{Q}}
    \quad\text{and}\quad
    P_{j+1}=\Theta_{j+1}^{\top}\begin{bmatrix}
      I_n & 0 \\
      0 &  0
    \end{bmatrix}\Theta_{j+1}
  \end{equation*}
\end{itemize}
Step 4.
\begin{itemize}
\item[] Set $j=j+1$ and goto Step 2.
\end{itemize}
\end{minipage}}\\

Algorithm 2 is almost identical to Algorithm 1 on the Gra{\ss}mannian,
except for the simpler Sylvester equation that is indeed a Lyapunov
equation here. Again, we immediately get the following result.

\begin{theorem}
  Algorithm 2 converges locally quadratically to any nondegenerate
  critical point of the Rayleigh quotient function on the Lagrange
  Gra{\ss}mannian. 
\end{theorem}

\subsubsection{Invariant subspace computation}
\label{sec:invariant-subspace-computation}  
We now turn to the more complicated task of solving the optimization problem
of the nonlinear trace function $\tr ((I-P)APA^{\top})$ over the Gra{\ss}mann
manifold $\Gr_{m,n}$. Here $A$ denotes an arbitrary real $n\times n$ matrix. This is
interesting as it leads to a locally quadratically convergent algorithm by
solving only linear matrix equations. 
Our method requires only orthogonal matrix
calculations and a linear matrix solver.  
We omit the straightforward
calculations that allow one to compute the Newton step 
in terms of the linear matrix equation appearing in Step 3 of the algorithm.

\hspace{-5mm}\fbox{\begin{minipage}[t]{125mm}
    \begin{center}
      \textbf{Algorithm 3: Invariant subspace function on the Gra{\ss}mannian}
    \end{center}
\noindent Step 1.
\begin{itemize}
\item[] Pick an orthogonal matrix $\Theta_{0}\in\SO_{n}$ corresponding
  to 
  \begin{equation*}
    P_{0}=\Theta_{0}^{\top}\begin{bmatrix}
      I_m & 0 \\
      0 &  0
    \end{bmatrix}\Theta_{0}\in\Gr_{m,n},
  \end{equation*}
  and set $j=0$.
\end{itemize}
Step 2.
\begin{itemize}
\item[] Compute
  \begin{equation*}
    \begin{bmatrix}A_{11} & A_{12} \\ A_{21} & A_{22}\end{bmatrix}=
  \Theta_{j}A\Theta_{j}^{\top}.
\end{equation*}
\end{itemize}
Step 3.
\begin{itemize}
\item[] Solve the linear matrix equation
  \begin{equation*}
   \begin{split}
   & A_{11}(A_{11}^{\top}Z_{j}-Z_{j}A_{22}^{\top})-
(A_{11}^{\top}Z_{j}-Z_{j}A_{22}^{\top})A_{22}\\
&-A_{21}^{\top}(Z_{j}^{\top}A_{12}+A_{21}Z_{j})-(A_{12}Z_{j}^{\top}+Z_{j}A_{21})A_{21}^{\top}
=A_{21}^{\top}A_{22}-A_{11}A_{21}^{\top}.
\end{split}  
\end{equation*}
for $Z_{j}\in\R^{m\times(n-m)}$.
\end{itemize}
Step 4.
\begin{itemize}
\item[] Compute
  \begin{equation*}
    \Theta_{j+1}^{\top}=\Theta_{j}^{\top}
    \begin{bmatrix}
      I_m & Z_{j} \\
      -Z_{j}^{\top} &  I_{n-m}
    \end{bmatrix}_{\mathrm{Q}}
    \quad\text{and}\quad
    P_{j+1}=\Theta_{j+1}^{\top}\begin{bmatrix}
      I_m & 0 \\
      0 &  0
    \end{bmatrix}\Theta_{j+1}
  \end{equation*}
\end{itemize}
Step 4.
\begin{itemize}
\item[] Set $j=j+1$ and goto Step 2.
\end{itemize}
\end{minipage}}\\

We do not address the interesting but complicated issue how to solve the above
linear matrix equation. Obviously one can always rewrite it as a 
linear equation on $\R^{m(n-m)}$  and then solve this, using matrix
Kronecker products and $\vecc$-operations, by any linear equation
solver. An alternative approach is to rewrite the equation in recursive
form as   

\begin{equation}
   \begin{split}
  A_{11}X_{j}-X_{j}A_{22}
&=A_{21}^{\top}(Z_{j-1}^{\top}A_{12}+A_{21}Z_{j-1})\\
&\quad+(A_{12}Z_{j-1}^{\top}+Z_{j-1}A_{21})A_{21}^{\top}-A_{21}^{\top}A_{22}+A_{11}A_{21}^{\top}\\
A_{11}^{\top}Z_{j}-Z_{j}A_{22}^{\top}&=X_{j}
\end{split}  
\end{equation}
starting from e.g. $Z_{0}=0$. This system of linear equations is uniquely
solvable if and only if the block matrices $A_{11}$ and $A_{22}$ have disjoint
spectra. 
Once again, we immediately get the following result. Recall, that an invariant
subspace $V$ of a linear operator $A:X\rightarrow X$ is called stable, if the restriction 
$A|_{V}$ and corestriction operators $A|_{X/V}$, respectively, have disjoint spectra.

\begin{theorem}
  Algorithm 3 converges locally quadratically to projectors onto stable invariant
  subspaces of $A$.
\end{theorem}

\section{Conclusions}
\label{sec:conclusions}
We presented a new differential geometric approach to Newton algorithms on a
Gra{\ss}mann manifold. Both the classical Gra{\ss}mannian as well as the
Lagrange Gra{\ss}mannian are considered. The proposed Newton
algorithms depend on the choice of a pair of local coordinate systems having
equal derivatives at the base points. Using coordinate charts defined by 
the Riemannian normal coordinates
and $QR$--factorizations, respectively, leads to an efficiently
implementable algorithm. Using the proposed method, new algorithms for symmetric 
eigenspace computations and non-symmetric invariant subspace
computations are presented that have potential for considerable computational
advantages, compared with previously proposed methods.      

\bibliographystyle{plain}

\end{document}